\documentclass{article}

\usepackage{amssymb}
\usepackage{graphicx}
\newcommand{\be}{\begin{equation}}
\newcommand{\ee}{\end{equation}}
\newcommand{\bea}{\begin{eqnarray}}
\newcommand{\eea}{\end{eqnarray}}

\begin{document}

\title{Open problems in mathematical physics}

\author{Alan A.~Coley}

%-----------------------------------------------------------------

%%\author{Alan~Coley}
%\footnote{aac@mathstat.dal.ca}
%%\vspace{5pt}
%%{\it {: {\small {aac@mathstat.dal.ca}}}
%%, Halifax, Canada.
%\pacs{PACS}%%\pacs{98.80.Es, 98.80.Cq}

\maketitle

\begin{abstract}

We present 
a list of open questions in mathematical physics. After a historical introduction,
a number of problems in a variety of different fields are discussed, with the intention
of giving an overall impression of the current status of mathematical physics, 
particularly in the topical fields of  
classical general relativity,  cosmology and the quantum realm.
This list is motivated by the recent article proposing
42 fundamental questions (in physics) which must be answered on the road to full enlightenment \cite{AL42}.
But paraphrasing a famous quote by the British football manager Bill Shankly, in response to the question of
whether mathematics can answer
the Ultimate Question of Life, the Universe, and Everything,
mathematics is, of course, much more important than that.

\end{abstract}

\vspace{30pt}

Alan Coley:
Department of Mathematics and Statistics, Dalhousie University,
Halifax, Nova Scotia, B3H 4R2, Canada [aac@mathstat.dal.ca].

%%============

%%one draft long--- pages with details re GR. book?where extended GR? liv. rev, arxiv, cqg
%%%Problems in GR--some in cosm/QGconsistent notation: GR throughout?

\newpage

\section{Mathematical Problems}

There are essentially two branches of mathematics, 
which in the broadest sense can be referred to as
pure mathematics and applied mathematics
(but there are actually three types of mathematicians; those that can count and those that cannot!). The 
actual mathematics (the problems, techniques and rigour) used in both are exactly the
same, but perhaps 
pure mathematicians and applied mathematicians are motivated differently.
Pure mathematics is concerned with mathematics for its own sake,
and an important criterion for assessing a worthy problem
is whether it leads to new developements in mathematics
[inwardly directed]. 
Applied mathematics is also (and perhaps primarily) concerned with establishing facts of real world interest
[outwardly directed]. For a more philosophical discussion on the nature of mathematics see, for example, the preface to
\cite{philosophy} and references within.

Noted probems in mathematics have always been important and are part of the
mathematical culture, both as recreation and as tests of acumen.
Unlike physics, 
where problems are dictated by necessity and practicalities, problems in mathematics,
particularly on the more pure side, have a life of their own and
the opinions of central characters have always been very 
important and played an elevated and pivotal role. Hence the importance
attached to problems espoused by famous mathematicians. 

In the sixteenth century, and according to the custom of the time, mathematical challenges, 
a type of intellectual duel
and a way of showing ones mathematical chops and gaining respect,  were often made. 
In 1530, there was a famous
contest between Niccol Tartaglia and 
Antonio Fiore (a student of
Scipione del Ferro) on solving 
cubic equations.  Each contestant proposed a number of problems for his rival
to solve, and whoever solved the most problems would receive 
all of the money put up by the two contestants. Since Tartaglia had worked out a general method
for solving cubic equations,
he won the contest.
Later, Tartaglia 
revealed his secret method to Gerolamo Cardano  
(which later led to questions of priority between 
Ferro and Tartaglia) when
Cardano published a book on cubic equations.
This subsequently led to a
challenge by Tartaglia, which was eventually accepted by Cardano's student Lodovico
Ferrari.  Ferrari beat Tartaglia in the
challenge, and Tartaglia lost both his prestige and income
\cite{Katz}.

In 1696 John Bernoulli published a challenging problem: {\em{To find the curve connecting two points, at different
heights and not on the same vertical line, along which a body acted upon only by gravity will fall in the shortest
time}} (the curve which solves this problem is called the `brachistochrone').  Gottfried Wilhelm Leibniz and Bernoulli were confident
that only a person who knew calculus could solve this problem (and it was rumoured that
this problem was set, in part, to determine what Isaac Newton knew on this topic since he had not
published his results yet).  Within one day of receiving the challange,  Newton sent in his solution.  When Bernoulli announced the winners of his contest, he named Leibniz
and l'Hopital (Leibniz's student) and one anonymous winner.  Bernoulli recognized the anonymous winner in
public with the phrase: ``we know the lion by his claw".

The twenty-three problems published by the mathematician David
Hilbert in 1900 \cite{Hilbert} are probably the most famous problems in mathematics. 
All of the problems were unsolved at the time of publication.
Several of them have been very influential in the development of mathematics.
Mathematicians and mathematical organizations have since announced several lists of problems, but
these have not had the same influence as Hilbert's original problems.
At the end of the millennium, which was also the centennial of Hilbert's publication of his problems, 
several mathematicians accepted the challenge to formulate
``a new set of Hilbert problems''.  
Most notable are Steven Smale's eighteen
problems, but to date these have not garnered very much popular attention.
Perhaps the twenty-first century analogue of Hilbert's problems 
is the list of seven Millennium Prize Problems chosen in the year 2000 by the Clay Mathematics Institute.

\newpage

\subsection{Hilbert's problems}

Hilbert's twenty-three problems in mathematics were published by  David Hilbert
in 1900 \cite{Hilbert}, and ranged over a number of topics in contemporary mathematics of the time.  Some of these problems were stated precisely
enough to enable a clear answer, while for others a solution to an accepted interpretation might have been possible
but closely related unsolved problems exist. And some of Hilbert's problems were not formulated precisely enough in
themselves, but were suggestive for more modern problems.  At the time of publication the problems were all
unsolved.  Several of them were very influential for twentieth century mathematics; for example, the 11th and the
16th problems (H11 and H16 -- see the text below and the Appendix where {\em{all}} of the problems
referred to are stated) have given rise to the flourishing mathematical subdisciplines of quadratic forms and
real algebraic curves.  A number of problems have given rise to solutions that have garnered great acclaim
including, for example, H1 and H10.  And many aspects of these problems are still of great interest today.

There are two problems that are not only unresolved but may, in fact, not be resolvable by modern standards.  
For example, H6 concerns the axiomatization of physics and H4 concerns the foundations of geometry. H4 
is generally thought to be too vague to enable a definitive answer,
and there is no clear mathematical consensus on the possible relevence 
of Godel's second incompleteness theorem (which gives a precise sense in which such a finitistic proof of the consistency of
arithmetic is unprovable). In addition,
Hilbert originally included a
``24th problem" (in proof theory, on a criterion for simplicity and general methods), but H24
was withdrawn from the 
list since it was regarded as  being too vague to ever 
be described as solved.

Noteworthy for its appearance on the list of Hilbert problems, and Smale's list and the list of Millennium Prize
Problems,  is {\em{the Riemann hypothesis}} (H8), which   asserts that
{\em{all nontrivial zeros of the analytical continuation of the Riemann zeta function have
a real part of 1/2.}}
A proof or disproof of this would have far-reaching implications in number theory.  
H8 is still considered to be an important open
problem, and has led to  other important
prime number problems, including Goldbach's conjecture and the twin prime conjecture, both of which remain unsolved.
However, even this famous hypothesis in pure mathematics is related to the energy eigenvalues of
distributions of random matrices, which is important in nuclear physics and quantum chaos \cite{Bender}.

\subsubsection{Summary and status of  Hilbert's problems}

Of the clearly formulated Hilbert problems, problems H3, H7, H10, H11, H13, H14, H17, H19, H20 and H21 
have a resolution that is generally accepted by consensus. On the other hand, problems H1, H2, H5, H9, H15, H18 and H22 have solutions that have been partially accepted, although there is some controversy as to whether the problems have been adequately resolved.

That leaves H8 (the Riemann hypothesis), H12 and H16 as unresolved.
H6 might be considered as a problem in physics rather than in mathematics. 
And H4 and H23 are too vague to ever 
be described as solved.

The 4 unsolved problems are \cite{Hilbert}:

\begin{itemize}

\item H6 Mathematical treatment of the axioms of physics.

\item H8 The Riemann hypothesis.

\item H12 Extend the Kronecker-Weber theorem on abelian extensions of 
the rational numbers to any base number field.

\item H16 Describe relative positions of ovals originating from a real algebraic curve and as limit cycles of a
polynomial vector field on the plane.  

\end{itemize}

The other Hilbert problems are listed in the Appendix. The majority of these problems are in pure mathematics;
only H19-H23 are of direct interest to physicists.
The Riemann hypothesis (H8), and H12 and H16 are problems in pure mathematics in the areas of number theory and algebra 
(and H16 is unresolved even for algebraic curves of degree 8).

H6 concerns the axiomatization of physics.
In particular, Hilbert proposed the following two specific problems:  (i) the axiomatic
treatment of probability with limit theorems for the foundation of statistical physics and (ii) the rigorous theory of
limiting processes ``which lead from the atomistic view to the laws of motion of continua."  Kolmogorov's 
axiomatics \cite{Kolmogorov}
is now accepted as standard and there has been some success regarding (ii) \cite{H6}. This is indeed a problem within mathematical physics, although
it is perhaps not necessarily regarded as being of prime importance 
in contemporary physics.

\newpage

\subsection{Smale's problems}

Steven Smale proposed  a list of eighteen unsolved problems in mathematics in 1998  \cite{Smale}, 
inspired by  Hilbert's original list of problems and at the behest of Vladimir Arnold.
Smale's problems  S1  and S13 are Hilbert's eighth (Riemann hypothesis)
and sixteenth (H8 and H16) problems, respectively, which remain unsolved.

The {\em{Poincare conjecture}} (S2), which asserts that
{\em{in three dimensions a sphere is characterized by the fact that it is the only closed and
simply-connected surface}},  
was proved by Grigori Perelman in 2003 using Ricci flows \cite{Perelman}. 
This problem is central to
the more general problem of classifying all 3-manifolds, and has many 
applications in modern theoretical physics.

There are 9 remaining  unsolved problems:

\begin{itemize} 
\item S3 Does P = NP?

\item S4 Shub-Smale conjecture on the integer zeros of a polynomial of one variable.

\item S5 Height bounds for Diophantine curves

\item S8 Extend the mathematical model of general equilibrium theory to include price adjustments.

\item S9 The linear programming problem:  find a strongly-polynomial time algorithm which decides whether, for given a matrix A 
(in $R^{m\times n}$)  
and b (in $R^{m}$),  there exists an x (in $R^{n}$) with $Ax \geq b$.

\item S10 Pugh's closing lemma (higher order of smoothness) 

\item S15 Do the Navier-Stokes equations in $R^3$ always have a unique smooth solution that extends for all time?  

\item S16 Jacobian conjecture  

\item S18 Limits of intelligence.

\end{itemize}

The famous problem {\em{Does P = NP?}}  (S3) is
{\em{whether or not, for all problems for which an algorithm can verify a given solution 
in polynomial time (termed a non-deterministic polynomial time or NP problem), an algorithm can also find that solution quickly
(a polynomial time or P problem)}}; that is,  whether all problems in NP
are also in P.  This is generally considered to be one of the most important open questions in mathematics and theoretical
computer science and it has far-reaching consequences to other problems in mathematics, and in biology, philosophy
and cryptography.  A common example of a P versus NP problem is the so-called
travelling salesman problem
(which asks the following question: {\em{Given a list of cities and the distances between each pair of cities, what is the shortest possible route that visits each city exactly once and returns to the origin city?}}) It is an NP-hard problem in combinatorial optimization, important in operations research and theoretical computer science.
Most mathematicians and computer scientists expect that the answer is that it is not true (i.e., $ P \neq NP$).
This problem also appears in the 
Millennium Prize list.

The problem S8 is in financial mathematics, which might be regarded as within the purview of theoretical physics.
Gjerstad
\cite{Gjerstad} has extended the deterministic model of price adjustment to a stochastic model and shown 
that when the stochastic
model is linearized around the equilibrium the result is the autoregressive price adjustment model used in applied
econometrics.  In tests it was found that the model
performs well with price adjustment data from
a general equilibrium experiment with two commodities.

Problems S4, S5, S9, S10 and S16 are problems in pure mathematics.
Smale also listed three additional problems in pure mathematics: the
Mean value problem, the question of whether 
the three-sphere is a minimal set, and
whether an Anosov diffeomorphism of a compact manifold topologically is the same as the Lie group model of John Franks?
The solved problems are listed in the Appendix. Unlike the Hilbert problems, many of these problems have practical applications and 
are of relevence in physics. For example, an
alternative formulation of S7 is the
Thompson Problem of the distribution of equal point charges on a unit sphere governed by the electrostatic Coulomb law. Problem S18 is concerned with the fundamental problems of intelligence and learning, both from 
the human and machine side.

The Navier-Stokes equations describe the motion of fluids. 
{\em{The problem is essentially to make progress towards a well-defined mathematical theory}} that will give insight into these equations.
Therefore S15 is truely a problem in mathematical physics and has imporant applications in many branches of theoretical physics
including engineering and oceanography, and even astrophysics.

Solutions of the compressible Euler equations typically develop singularities
(that is, discontinuities of the basic fluid variables), in a finite time \cite{Sideris}.
The proofs of the development of singularities are often by contradiction
and consequently do not give detailed information on what occurs when the smooth solutions
break down. The formation of shock
waves are possible, and it is known that  in some cases solutions can be physically extended
beyond the time of shock formation. The extended
solutions only satisfy the equations in the ``weak sense''. For the classical Euler
equations there is a well-known theorem on the global existence of weak solutions
in one (space) dimension \cite{Glimm}, and a  one dimensional
class of weak solutions  has recently been 
found in which both existence and uniqueness hold \cite{Bressan}. In higher (space) dimensions there are no general global existence
theorems known. The question of which quantities must blow up
when a singularity forms in higher dimensions has been partially addressed
for classical hydrodynamics \cite{Chemin}.
A smooth solution of the classical Euler equations 
has been proven to exist for all time when the initial data are small and the fluid is
initially flowing outwards  uniformly \cite{Grassin}.

\newpage

\subsection{Millennium Prize problems}

The Millennium Prize Problems are seven problems in mathematics that were proposed by 
the Clay Mathematics Institute
in 2000 \cite{Millennium}, with a \$1 million US prize being awarded by the Institute to the
discoverer(s) of a
correct solution to any of the problems.
At present, the only Millennium Prize problem to have been solved is the Poincare conjecture \cite{Perelman}.
In addition to the 
Poincare conjecture, 
three other  problems, namely the Riemann hypothesis (H8), P versus NP (S3), and the existence and smoothness
of the Navier-Stokes equations (S15), are also on Smale's list.

There are 3 remaining unsolved problems \cite{Millennium}:

\begin{itemize}

\item M1 The
Hodge conjecture that for projective algebraic varieties, Hodge cycles are rational linear combinations of algebraic cycles.

\item M2
Yang-Mills existence and mass gap.

\item M3 The
Birch and Swinnerton-Dyer conjecture.

\end{itemize}

Problem M2 aims to establish {\em{the existence of the quantum Yang-Mills
theory and a mass gap}} rigorously, and is truely a problem in mathematical physics.
Classical Yang-Mills theory \cite{YangMills} is a generalization (or
analogue) of Maxwell's theory of
electromagnetism in which the chromo-electromagnetic field itself carries charges.  As a classical field theory,
it's solutions
propagate at the speed of light and so its quantum version describes massless  gluons.  
The so-called mass
gap is the problem that
color confinement only allows bound states of gluons, which form massive particles. 
The asymptotic freedom of confinement also makes it possible that a quantum Yang-Mills theory exists
without restriction to low energy scales.

Many important mathematical
questions remain unsolved, including 
stability theorems and the proof of existence of
Yang-Mills fields by methods of partial differential equations.
More contemporary 
questions are to obtain solutions of the Yang-Mills equations on a Riemannian (or Lorentz)
manifold.  The Yang-Mills equations in general relativity will be discussed later.

The Birch and Swinnerton-Dyer
conjecture M3 asserts that
{\em{that there is a simple way to tell whether the equations defining elliptic curves have a finite or infinite 
number of rational solutions}}.  This is a special case of
Hilbert's
tenth problem, in which it has been proven that there is no way to decide whether
a given equation in the more general case even has any solutions.

\newpage
\section{Mathematical Physics}

Not all mathematical problems are necessarily of interest to a physicist. Similarly,
not all problems in physics are of a mathematical character. 
For example, there are many lists of problems in physics, including problems in 
high-energy physics/particle physics,
astronomy and astrophysics,
nuclear physics,
atomic, molecular and optical physics,
condensed matter physics, and
biophysics \cite{Listphysics,Baez}. But these cannot all be regarded as problems within mathematics.
Most problems of a mathematical nature are restricted to fundamental physics and particularly
theoretical physics (and especially in
theories such as general relativity (GR) and quantum gravity (QG)). It is perhaps illuminating to recall the quote by Werner von Braun who said that
``Basic research is what I am doing when I don't know what I am doing''.

Five of the most important and interestingly unsolved problems in theoretical physics in the 
quantum regime (in the small) and in cosmology (in the large) are commonly agreed to be the following
(see, for example, \cite{gonitsora}):

\begin{itemize}

\item Ph1 The Problem of Quantum Gravity.

\item 
Ph2  The Foundational Problems of Quantum Mechanics.
\item 
Ph3  The Unification of Particles and Forces.
\item 
Ph4 The Tuning Problem.
\item 
Ph5 The Problem of Cosmological Mysteries.

\end{itemize}

We shall be interested in problems which we shall refer to as problems in mathematical physics, which we shall
define to mean problems that are 
well-formulated (i.e., well-posed) mathematical problems, which are of interest to physicists.
Many such problems involve 
systems of partial differential equations, which are of central importance in theoretical physics.

In general, problems in mathematical physics will not include problems 
where the
basic underlying physics is not understood (such as, for example
the quantization of gravity), and although it is clear that their
solution will inevitably involve a lot of mathematics (and perhaps even lead to new areas of mathematics),
an explicit well-posed mathematical problem cannot be formulated. Nor do they
include problems in pure mathematics where there is no clear physical application
(e.g., the Riemann hypothesis). There are also questions in 
computational mathematics, and it is also debatable whether such problems qualify as
problems  in  mathematical physics.
The meaning of 
problems in mathematical physics is nicely
illustrated by the set of 15 open problems proposed by mathematical physicist Barry Simon 
\cite{simon}, which we shall discuss a little later.

\newpage

This paper is  motivated, in part,  by the recent article entitled
{\em{Life, the Universe, and everything: 42 fundamental
questions}} 
(referred to hereafter as {\bf{AL42}} \cite{AL42}: the
actual list of questions is given in the appendix), which itself was inspired by
{\em{The Hitchhiker's Guide to the Galaxy}}, by Douglas Adams.
There are many questions in theoretical physics discussed in {\bf{AL42}}, some of which are
of relevance here and will be discussed in more detail later, including
the cosmological constant problem (AL2.1), the dark energy problem (AL2.2),
the regularization of quantum gravity (AL2.3), black hole entropy and thermodynamics (AL2.4),
black hole information processing (AL2.5),
supersymmetry and the hierarchy problems (AL3.3), and
higher dimensions and the geometry and topology of internal space (AL5.1).

In this paper, I shall present a number of what I consider to be problems in
mathematical physics, primarily in the current areas of theoretical and fundamental physics.
Classical GR  remains healthy and vigorous, in
part, due to a frequent injection of fertile mathematical ideas (such as those of
Hawking and Penrose and, more recently, of Schoen-Yau and Witten). By any
reasonable definition of the term, it is clear that much of classical GR
is ``mathematical physics".
GR problems have typically been under-represented in lists of problems in mathematical physics (e.g., see  \cite{simon}), perhaps due to their advanced technical nature. Obviously any such list is subjective, and 
classical GR may well be over-represented here, but I feel at liberty to comprehensively
discuss problems in GR (artistic licence?) and
to present some of my own
personal favorites ({\bf{PF}}) (perhaps to justify my own research interests?).

After the current more introductory and historical section,
I shall discuss in more detail classical GR first, and then return to quantum theory and cosmology
(and specifically discuss the 5 physics problems above) in the ensuing sections.
It is the technical problems that are of interest to mathematicians. Often physicists
are perhaps not as interested in the technical aspects of the problem, but more 
in the context and the consequences of the results. Hence, although I shall attempt
to state the problems relatively rigorously, as is appropriate for mathematicians,
I shall endeavor to select less technical questions, or at least describe them in as heuristic manner as possible,
which may well be of more interest to 
physicists.

This article is written primarily for a
readership with some background in mathematics and physics.
However, regardless of 
background, the intention here is not for readers to
understand each and every problem, but rather to get an overall 
impression
of the open questions in the various fields. In particular, one aim is 
to outline which
areas are currently exciting with  unsolved problems whose potential 
solution might have a huge impact on the field,
and consequently motivate readers (and especially young physicists) to 
possibly get more involved in research.
Obviously this article takes it for granted that mathematics is necessarily 
the language of physics
(that is, the so-called {\em{unreasonable effectiveness of mathematics 
in the natural sciences}} 
\cite{Wigner});
however, the philosophical reasons for this it is beyond the current 
discussion.

\newpage

\subsection{More on lists}

There are many lists of 
unsolved problems in mathematics (see, for example, \cite{Listofunsolved} \cite{AbeTanaka}).
These include many problems in applied mathematics (and hence mathematical
physics), some of which have been discussed above (the regularity of the  Navier-Stokes and
Yang-Mills equations, and  problems on turbulence).
In particular, there are questions on stability (e.g., for what classes of 
ordinary differential equations, describing dynamical systems, does the Lyapunov second method formulated in the classical and canonically generalized forms define the necessary and sufficient conditions for the asymptotic stability of motion?), questions in  ergodic theory
(e.g., the Furstenberg conjecture), on  actions in higher-rank groups
(e.g., the Margulis conjecture), the question of whether
the Mandelbrot set is locally connected, and problems in Hamiltonian flows (e.g., the
Weinstein conjecture: Does a regular compact contact type level set of a Hamiltonian on a symplectic manifold carry at least one periodic orbit of the Hamiltonian flow?)

In particular, very recently the 
DARPA Mathematical Challenges were proposed \cite{DARPA}, which are very heavy in applied mathematics and theoretical physics. They involve not only problems in
classical fluid dynamics and the Navier-Stokes Equation (and their use in the quantitative understanding of shock waves, turbulence, and solitons), but also problems in which new methods are needed to tackle complex fluids (such as foams, suspensions, gels, and liquid crystals), and  the  Langlands program (see below).
In addition, a number of DARPA challenges involve
traditional problems in pure mathematics, such as the Riemann hypothesis 
(number theory), the Hodge conjecture (in algebraic geometry),
and in convex optimization (e.g., whether linear algebra be replaced by algebraic geometry in a systematic way).
They also include
the  physical consequences of Perelman's proof of Thurston's geometrization theorem and 
the implications for spacetime  and cosmology of the Poincare conjecture in four dimensions.

Also a number of more speculative problems were proposed in an attempt to apply mathematics to new areas of interest, including
the mathematics of the brain, the dynamics of networks, stochasticity in nature, problems in
theoretical  biology  and biological quantum field theory
(e.g., what are the fundamental laws of biology, can  Shannon's information theory  
be applied to virus evolution, the geometry of genome space, what are the symmetries and action principles for biology) and the mathematics of quantum computing (algorithms and entanglement)
including optimal nanostructures, and problems in theoretical
computation in many dimensions.
One of the most important advances in the last few years has been the use
of theoretical computing and  neural networks to attempt to solve all kinds of previously untractable
problems.

There are also a number of interesting questions, some of which are discussed in
{\bf{AL42}}, which might be considered to be more metaphysics than physics,
and certainly outside the realm of mathematical physics (although they may be
addressed by scientists, and indeed mathematicians, in the future).
These include the study of the multiverse and the anthropic principle,
and emergent phenomena such as life and
consciousness (the puzzle of the possible role of human consciousness in resolving questions in
quantum physics is discussed by \cite{Penrosebooks}). For example, in 
Section 7 of {\bf{AL42}}
the ultimate nature of reality, the reality of human experience,
conscious minds and questions on the origin of
complex life are broached. To this list, questions of 
ethics and even religion might be added. The 
potential for breakthroughs in theoretical, computational, experimental, and
observational techniques are also discussed in {\bf{AL42}}. Although
such topics are outside the purview of the current article, that is not to
say that mathematics might not be useful in their consideration.

\newpage

\subsection{Mathematical physicists}

In mathematics, the Langlands program \cite{Langlands} constitutes a number of conjectures that relate Galois groups in algebraic number theory to automorphic forms and representation theory of algebraic groups over local fields. DARPA proposed two challenges: (geometric Langlands and quantum physics)
how does the Langlands program
explain the fundamental symmetries of physics (and vice versa), and (arithmetic Langlands, topology, and geometry)
what role does homotopy theory play in the classical, geometric, and quantum Langlands programs.

It has been thought for a long time
that the Langlands duality ought to be related to various dualities observed in
quantum field theory and string theory. The so-called Langlands dual group \cite{Langlands}, which is
essential in the formulation of the Langlands correspondence,  plays an important role
in the study of S-dualities in physics and was introduced by  physicists 
in the framework of four-dimensional gauge theory \cite{Goddard}.
Witten recently showed that
Langlands duality is closely related to the S-duality of quantum field theory,
which opens up exciting possibilities for both subjects \cite{KapustinWitten}.
Indeed, the connections between the Langlands
program and two-dimensional conformal field theory give important insights into the physical implications of the Langlands duality.

Edward Witten is a theoretical physicist working in string theory, quantum gravity, supersymmetric quantum field theories, and other areas of mathematical physics.
In addition to his contributions to physics, Witten's work has also significantly impacted pure mathematics. In 1990 he became the first (and so far only) physicist to be awarded a Fields Medal by the International Mathematical Union. 
The Fields Medal is regarded as the highest honour a mathematician can receive and, together with the Abel Prize, has often been viewed as the ``Nobel Prize" for mathematics. 

In a written address to the International Mathematical Union, Michael Atiyah said of Witten \cite{Atiyah}:
``Although he is definitely a physicist his command of mathematics is rivaled by few mathematicians, and his ability to interpret physical ideas in mathematical form is quite unique''.
As an example of Witten's work in pure mathematics, Atiyah cited his application of techniques from quantum field theory to the mathematical subject of low-dimensional topology. In particular, Witten realized that Chern–-Simons theory in physics could provide a framework for understanding the mathematical theory of knots and 3-manifolds \cite{Witten89}.
Witten was also awarded the Fields Medal, in part, for his proof in 1981 of the positive energy theorem in general relativity \cite{Witten81}.
 
There are also many mathematicians who have greatly influenced physics. These include Roger Penrose and Steven Hawking (whose contributions will be discussed later).
Michael  Atiyah is a mathematician specializing in geometry, and was awarded the Fields Medal in 1966. He helped to lay the foundations for topological K-theory, an important tool in algebraic topology. 
The Atiyah--Singer index theorem \cite{AtiyahSinger} 
(in which the index is computed by topological means) is widely used in counting the number of independent solutions to differential equations. The index theorem provides a link between geometry and topology 
and has many applications in theoretical physics.
Some of his more recent theoretical physics inspired work, and particularly that on instantons and monopoles,  is responsible for some subtle corrections in quantum field theory.

Simon Donaldson, one of Atiyah's students, 
is known for his work on the topology of smooth (differentiable) four-dimensional manifolds and the Donaldson (instanton) invariant (among other things). 
Donaldson's work is on the application of mathematical analysis (and especially that of elliptic partial differential equations) to problems in the geometry of 4-manifolds, complex differential geometry and symplectic geometry \cite{Kronheimer}.
He has used ideas from physics to solve mathematical
problems, and investigated problems in mathematics which have physical applications (e.g., 
an application of gauge theory to four-dimensional topology \cite{Donaldson}). Recently,
Donaldson's work has included a problem in complex differential geometry regarding a conjectured relationship between the stability conditions for smooth projective varieties and the existence 
of Kahler-Einstein metrics with constant scalar curvature
\cite{ChenDonaldsonSun}, which is of interest in string theory. String theory is often described as a topic within mathematics rather than in physics (in much the same way GR was fifty years ago).

The mathematician
Shing-Tung Yau was awarded the Fields Medal in 1982. 
Yau's work is mainly in differential geometry, especially in geometric analysis. He has been active and very influential at the interface between geometry and theoretical physics (see later). 
Together with Schoen, Yau used variational methods to prove
the positive energy theorem in GR, which
asserts that (under appropriate assumptions) the total energy of a gravitating system is always positive and can 
vanish only when the geometry is that of flat Minkowski spacetime. It consequently establishes Minkowski space as a stable ground state of the gravitational field. 
As mentioned above,
Witten's later simpler (re)proof  \cite{Witten81} used ideas from supergravity theory.
Yau also
proved the Calabi conjecture which allows physicists to demonstrate, utilizing Calabi–-Yau compactification, that string theory is a viable candidate for a unified theory of nature. Calabi-Yau manifolds are currently one of the standard tools for string theorists.

\newpage
\subsection{Simon's problems}

Problems in mathematical physics are 
well formulated mathematical questions of interest to physicists.
The meaning of problems in mathematical physics is nicely
illustrated by the set of 15 open problems proposed by mathematical physicist Barry Simon in 1984 \cite{simon},
who was awarded the American Mathematical Society's Steele Prize for Lifetime achievements in mathematics
in 2016.

I shall display and briefly discuss six of these problems below.
The first two questions are in fluid dynamics and have been alluded to earlier.
The sixth, cosmic censorship, will be discussed later.
The
remaining problems are displayed in the Appendix (the citations therein are  circa 1984  \cite{simon},
and there has subsequently been progress on these problems). Although many of these problems 
involve Schrodinger operators, Simon's own field of expertise, I believe that the problems do help give
a flavour of what problems constitute mathematical physics to a general physicist (for example,
one who is not necessarily 
an expert in GR, one of the fields to be discussed below).

\begin{itemize}

\item
BS1 Existence for Newtonian Gravitating Particles.
A: Prove that
the set of initial conditions  which fails to have global solutions is of
measure zero (some mathematicians believe that there may be
an open set of initial conditions leading to non-global solutions).
B: Existence of non-collisional singularities in the Newtonian
N-body problem.

\item
BS3 Develop a comprehensive theory of the long
time behavior of dynamical systems including a theory of the onset of, and of
fully developed, turbulence.

\item
BS2 Open Questions in Ergodic Theory. Particular problems include
A: Ergodicity of gases with soft cores, B: Approach to equilibrium, and C: Asymptotic Abelianness for the quantum Heisenberg dynamics.

\item
BS8 Formulation of the Renormalization Group and Proof of Universality.
A: Develop a mathematically
precise version of the renormalization transformations for $\nu$-dimensional
Ising-type systems.
B: In particular, show that the critical exponents in the
three dimensional Ising models with nearest neighbor coupling but different
bond strengths in the three directions are independent of the ratios of these bond
strengths.

\item BS14 Quantum field theory remains a basic element of fundamental
physics and a continual source of inspiration to mathematicians.
A: give a precise mathematical construction of quantum chromodynamics,
the model of strong interaction
physics. B: Construct any non-trivial renormalizable
but not super--renormalizable quantum field theory. C: Prove that quantum electrodynamics is not a consistent
theory. D:  Prove that a non-trivial lattice cutoff theories theory does not
exist.

\item
BS15 Cosmic censorship.

\end{itemize}

Problem BS3 is very general and rather vague, and so
the first problem is to formulate the really significant questions. For
recent reviews of some of the more spectacular developments see \cite{Feigenbaum,Collet80}.
There has been considerable progress in
understanding the onset of turbulence (e.g., see \cite{Ruelle}),
but fully developed turbulence is far from being comprehensively understood.
Even the connection between turbulence and the Navier-Stokes equation is
not absolutely clear
\cite{Foias}.

Regarding BS2, the developers of statistical mechanics and thermodynamics, including Boltzmann and
Gibbs, realized that from a microscopic
point of view  bulk systems rapidly approach equilibrium states parametrized by a 
few macroscopic parameters. It was originally believed that it could be proven that the classical dynamics
on the constant energy manifolds of phase space is ergodic.
However, the Kolmogorov--Arnold--Moser (KAM) theorem \cite{KAM} is a result in dynamical systems about the persistence of quasiperiodic motions under small perturbations. An
important consequence of KAM is that many classical systems will not be
ergodic: there will be an invariant subset of phase space consisting of a union of
invariant tori of positive total measure.

Problem BS14 concerns the question of whether quantum
field theory really is a mathematical theory at all.
This question remains open for any nonlinear
quantum field theory in three-space plus one-time dimensions.
The basic difficulty in formulating the mathematical problem is the
singular nature of the nonlinear equations proposed.
Physicists eventually developed sets of ad hoc rules to cancel
the infinities in QFT
and to calculate observable effects. These rules of renormalization were remarkably
accurate in producing verifiable numbers in electrodynamics.

Fisher, Kadanoff and Wilson \cite{Wilson74}  developed  the ``renormalization group theory'' of critical phenomena
which, regarding question
BS8.B above,
is often claimed to ``explain'' universality
(rather than universality being assumed).
The basic idea of shifting
scales as one approaches a critical point via a nonlinear map of Hamiltonians
and obtaining information from the fixed points of that map has been applied in a
variety of situations \cite{Feigenbaum}.
In some studies, the nonlinear maps are on well defined
spaces and there has been considerable progress on a rigorous mathematical
analysis on the Feigenbaum
theory \cite{Collet81}. The original Wilson theory is on functions of infinitely many
variables and it is far from clear how to formulate the maps in a mathematically
precise way (let alone then analyze their fixed point structure);  indeed, there are
various no-go theorems \cite{Griffiths79} on how one might try to make a
precise formulation in 
lattice systems.

\newpage

\subsection{Yau, Penrose and Bartnik}

Analytical methods (and especially the theory of partial differential equations) used 
in the study of problems in differential geometry, and subjects
related to geometry such as topology and physics, were surveyed
in  \cite{Yau1982}.
There was a section  in \cite{Yau1982} with       
120 open questions by Yau himself (p669).
Most of these problems are technical and in 
differential geometry (and mostly Riemannian geometry), and
are old and well known (even in 1982; see
original references therein). Many of the problems are not
related to physics directly, and hence are not necessarily 
problems in mathematical physics. But some of the problems concern
the Dirac equation,
gravitational instantons, Kahler and Calabi manifolds and Gauss--Bonnet theory.
There were 2 problems in Yang-Mills theory
(Problems Y117 and Y118), 
and 5 problems in GR: problems Y115, Y116 \& Y119 concern the topology
of a geodesically complete Lorentzian 4-manifold of nonnegative Ricci
curvature which contains an absolutely maximizing
timelike geodesic (see later),  the topology
of a static stellar model, and the characterization
of asymptotical flatness of a manifold in terms of a suitable decay rate of the  curvature,
respectively. The problems
Y114 (cosmic censorship) and
Y120 (the definition
of total angular momentum) are also included in the list of open problems by Penrose
(RP12 \& RP10, respectively) in
the same book \cite{Yau1982}.

The fourteen unsolved problems in classical
GR presented by
Roger Penrose (p631 in \cite{Yau1982}; problems RP1 -- RP14 in the Appendix), 
represented the status of the subject
circa 1982. (An 
earlier list of 62 problems in GR was given by Wheeler
\cite{wheeler}). Many of them were technical questions concerning
definitions of null infinity, appropriate (conformal) properties, and conservation laws and physical quantities, necessary for the formulation of the important problems and conjectures that followed.
In particular, in 1982 it was known that spherically symmetrical collapse models lead to
a black hole horizon, but if the initial data is perturbed
away from spherical symmetry, a so-called naked singularity could arise (from which causal curves
can extend to external future infinity). But the belief
was that naked singularities will not arise ``generically", whence it is said that cosmic censorship
holds \cite{Penrose1969,Penrose1978}.

Problem RP11, which is related to problem RP4, is necessary for the
statement of the  cosmic censorship problem RP12, which was stated
somewhat vaguely; indeed,
it is a problem in itself to find a satisfactory mathematical formulation
of what is physically intended \cite{Penrose1978,Penrose1979} (such as, for example,
are ``generic" maximally extended
Ricci-flat spacetimes globally hyperbolic or necessarily have a Cauchy surface
\cite{Geroch}). 

With a suitable assumption of cosmic censorship, together with
some other reasonable physical assumptions, it is possible to derive a certain
sequence of inequalities \cite{Penrose1973}. 
Problem
RP13 concerns the Penrose inequality, which
generalizes RP6 and is related to RP7. The validity of these inequalities are sometimes regarded as giving some credence to
cosmic censorship.

There are many other problems involving black holes which have not
yet been solved, including RP14. In particular, there are many open problems generalizing vacuum results to results with matter.
Generally
results for the Einstein-Maxwell equations are similar to those for the
pure Einstein vacuum equations, and Einstein-Maxwell analogues exist
for the problems RP3, RP4, RP9 and RP11. However, the statement of  problem RP14
is not true in the presence of electromagnetic fields.

There is also a list of open problems in mathematical GR by 
Robert Bartnik \cite{Bartnik} (also see references within).
Theoretical GR had developed to such an extent that rigorous mathematical arguments 
have replaced
many of the formal calculations and heuristics of the past, which
will yield new insights for both mathematics and physics.
Many of the Bartnik  \cite{Bartnik}
problems are
technical and concern clarifications and motivations for important contemporary
problems, and many have been noted elsewhere in this paper.

The problems are on the topics of (i)
apparent horizons (RB1 - RB17), (ii)
initial data sets (RB8 - RB112), (iii)
uniqueness and rigidity theorems for static and stationary metrics (RB13 - RB17),
(iv) approximations (RB18 - RB25), (v)
maximal and prescribed mean curvature surfaces (RB26 - RB29), (vi)
causality and singularities (RB30 - RB34), (vii)
the initial value problem and cosmic censorship (RB35 - RB47), and (viii)
quasi-local mass (RB48 - RB53).

Regarding (iv), there has been a lot of work done on constructing metrics
which approximately satisfy the Einstein equations, primarily consisting of 
numerical computation, but also involving asymptotic
expansion, linearisation and matching techniques. As noted earlier, it is debatable
as to whether numerical problems are in the realm of mathematical physics. But problem
RB21 concerns a rigorously proof of the Newtonian  limit to the Einstein equations
and problem RB20 concerns the range of validity of post-Newtonian and post-Minkowskian asymptotic expansions.   Problem RB23 on whether test particles  follow spacetime
geodesics,  is a famous problem and includes an extensive investigation
of asymptotic expansions \cite{JEhlers}.

Problem RB32 in (vi) is the  ``Bartnik splitting
conjecture'': {\em{Let M be a ``cosmological spacetime'' satisfying the timelike
convergence condition: then either M is timelike geodesically incomplete or M
splits as $R \times M^3$ isometrically (and thus is static).}}
This is essentially problem Y115 in \cite{Yau1982}, which posed the question of establishing a Lorentzian
analogue of the Cheeger-Gromoll splitting theorem
of Riemannian geometry \cite{CG}. The concept of geodesic completeness in Lorentzian
geometry differs considerably from that of Riemannian geometry, and this
question was  concretely realized in the Bartnik splitting conjecture  RB32 \cite{Bartnik}.
In the case of a 4D vacuum (i.e., Ricci flat) globally hyperbolic,
spatially compact
spacetime, if M splits it is
necessarily flat and covered by $R \times T^3$,
and thus for a non-vacuum ``cosmological
spacetime"  the conjecture asserts that the 
spacetime either is singular or splits. 
The resolution of the basic
Lorentzian splitting  conjecture  as considered in RB32
was given in \cite{EJH}, and can be viewed as a (rigidity) singularity theorem 
since the exceptional possibility that spacetime splits  can be
ruled out as unphysical, and hence the spacetime has an inextendible time–like geodesic which ends after a finite proper time (i.e.,
it is time–like geodesically incomplete and hence singular).
The status of the Bartnik conjecture was discussed in \cite{Galloway96}, and
more general versions of the
conjecture and partial results were discussed in \cite{GallowayVega}.

Regarding (vii), Bartnik stated there are many versions of cosmic censorship, but that essentially the aim is to 
prove a theorem showing that singularities
satisfying certain conditions are not naked.
In addition, problem RB43 concerns the 2-body system in Einstein gravity, which Bartnik claimed is
probably the most embarassing indictment of our (lack of) understanding of the
Einstein equations (however, see  \cite{Focusissue} and the discussion later). The
problem in (viii) of defining the total energy of an isolated system was essentially solved
in \cite{ADM}, but the correct definition of the energy content of a bounded region in
spacetime is still not settled. Although a number of candidate definitions have been suggested, so far  
none of these verify all
the properties expected of a quasi-local mass.

There have been a number of reviews on the global existence problem in GR,
including those of  \cite{Rendall2002} and  \cite{LARS99} (also see references within). In these reviews there is an 
emphasis on very technical questions in differential
geometric and analytical  global properties of 1+3 dimensional spacetimes containing a compact
Cauchy surface  (and particularly
the vacuum case), but they draw attention to a number of open questions
in the field.

\newpage
\section{Open problems in General Relativity}

Mathematical questions about the general properties of solutions of Einstein's field equations of GR
are truely problems in  mathematical physics. Problems in GR are
not necessarily more important than other problems in theoretical physics, but they do often have a more well-formulated mathematical expression. They are also perhaps more difficult for a broad based physics 
audience to fully appreciate. Therefore, I will first review some mathematical background, which can be 
skipped by general readers.

In general,  a smooth (or sufficiently differentiable)
4-dimensional Lorentz manifold (M, g) is considered.
The Lorentzian metric, g, which defines the causal structure on M, is
required to satisfy the Einstein
field equations, which constitute a 
hyperbolic system of  quasilinear partial differential
equations which are, in general, coupled to other partial differential equations describing
the matter content of spacetime \cite{Rendall2002}.
Primarily the vacuum case (when g is Ricci
flat) is considered.
Physicists are then interested in the Cauchy problem in which
the unknowns in the resulting Einstein vacuum constraint equations,  consisting of a Riemannian metric
and a symmetric tensor defined on a three-dimensional manifold (and initial data for any matter fields
present), are the initial data for the remaining Einstein vacuum evolution
equations.

The Einstein equations
are invariant under a change of the coordinate
system (general covariance or gauge freedom), which
complicates the way they must be formulated in order to 
faciliate the study of their global properties \cite{LARS99}.
Although the Einstein vacuum equations are not hyperbolic in the usual sense
due to general covariance, 
the Einstein vacuum
equations in spacetime harmonic coordinates constitute a quasi–linear hyperbolic
system and therefore the Cauchy problem is well posed and
standard results imply local existence  \cite{CB69}. It is also possible to show that if the constraints and
gauge conditions are satisfied initially, they are preserved by the evolution.
For example, the global regularity and modified scattering 
for small and smooth initial data with suitable decay at infinity
for a coupled Wave-Klein-Gordon system  
(a simplified version of the full Einstein-Klein-Gordon system) in 3D was studied in \cite{Ionescu}.
Analogues of the results for the vacuum Einstein equations  
are known for the Einstein equations coupled to many different types of matter, including perfect fluids, 
gases satisfying kinetic theory, scalar
fields, Maxwell fields, Yang-Mills fields and various combinations of these. 

The general results for perfect fluids
only  apply in the restricted circumstances in which the energy density is uniformly bounded away from zero
(in the region of interest) \cite{Rendall2002}.
The existence of global solutions for models with more exotic matter, such as stringy matter,
has also been studied \cite{Narita}.

\newpage

\paragraph{Existence:}
The basic local existence theorem says that, given smooth (i. e., infinitely differentiable  $C^{\infty}$) data for the vacuum Einstein equations, there exists a smooth solution of the equations (on a finite time interval) which
gives rise to these data \cite{CB80}.
The standard global uniqueness theorem for the Einstein
equations asserts that the long term solution (maximal development \cite{CB69}) of any Cauchy
data is unique up to a diffeomorphism which fixes the initial hypersurface and that, in an
appropriate sense, the solution depends continuously on the initial data \cite{CB80}.

The local existence of solutions of the Einstein equations is 
understood quite well. However, the problem of proving general global existence theorems
for the Einstein equations is beyond the reach of current mathematics \cite{Rendall2002}. 
The usual method for solving the Einstein equations is the conformal
method \cite{CB80}, in which the so-called free data are chosen
and the constraints then reduce to four elliptic equations.
In the simplified constant mean curvature 
case these reduce further to a linear system of three elliptic equations, which decouple
from the remaining equation which essentially reduces to the  nonlinear, scalar Lichnerowicz equation.

The causal structure of a Lorentzian spacetime is conformally invariant.
Friedrich derived the compactified ``regular conformal field equations''
from the Einstein
equations, 
a first order symmetric hyperbolic system, which leads to
well posed evolution equations and hence 
small data global existence
results from the stability theorem for quasi–linear hyperbolic equations.
For example,  Friedrich \cite{Friedrich86} proved global existence to the future for ``small''
hyperboloidal initial data  (that is, data close to the standard data on a hyperboloid) in Minkowski
space.
It is still an open question what general conditions on initial data on an
asymptotically flat Cauchy surface give a Cauchy development with regular
conformal completion. Friedrich has developed an
approach to this problem in which the conformal structure at spatial
infinity is analyzed (see \cite{Dain} for references, and  \cite{Kroon}
which points out some new obstructions to regularity; also see the more recent articles
\cite{HansRingstrom17,Ringstrom2015} and references within).

Therefore, for the full 1+3 dimensional Einstein equations (without symmetries) the only global existence
results known are the theorem on nonlinear stability of Minkowski space \cite{ChristodoulouKlainerman90}, the 
semi–-global existence theorem  for the hyperboloidal initial value problem \cite{Friedrich86} and the semi-–global existence
theorem for spatially compact spacetimes with Cauchy surface of
hyperbolic type \cite{AnderssonMoncrief}, which are all small data
results. It has been shown that for analytic vacuum or electrovac spacetimes,
with an analytic Cauchy horizon which is assumed to be
ruled by closed null geodesics,  there exists a nontrivial
Killing field \cite{Isenberg85}. Theorems in the cases of special spacetimes with symmetries are briefly reviewed below. Since  spacetimes with Killing
fields are non–generic, this result may be viewed as supporting evidence for the
strong cosmic censorship (see below).

\newpage
\paragraph{Special cases:}
It is possible to solve the global existence
problem for the Einstein equations in special cases, such as for spacetimes with
symmetry  \cite{Rendall2002,LARS99}.
For example,
basic global existence theorems for spherically symmetric static solutions (which are
everywhere smooth) 
have been proved for perfect fluids and  collisionless matter (see \cite{Rendall2002} and references within).
The spacetime symmetry is defined by the number
and character of Killing vectors.
For example, consider spacetimes with an $r$–-dimensional Lie
algebra of space–like Killing fields. For each $r \leq 3$, there are some basic
results and conjectures on global existence and cosmic censorship \cite{LARS99}.                     
In the cases $r$ = 3 (Bianchi models;
see, for example, \cite{Chrusciel55}) and a special case of $r$ =
2 (polarized Gowdy models -- see Refs. below), the global behavior of the Einstein
equations is well understood.

For the general
$r$ = 2 case (local $U(1) \times U(1)$ $G_2$ symmetry), there are only partial
results on the global existence problem and the cosmic censorship problem
remains open \cite{LARS99}.
The first global existence result for Gowdy spacetimes
with topology $R \times T^3$ was proven in  \cite{Moncrief81}, and subsequently generalized for spacetimes on $S^3$
and $S^2 \times S^1$ in \cite{Chrusciel47} (a class of ``nongeneric'' metrics still remains to be studied). The first result concerning global constant mean curvature foliations in  vacuum Gowdy spacetimes
was proven in \cite{IsenbergMoncrief1982}.
The question of cosmic censorship for the Gowdy spacetimes may be studied
by analyzing the asymptotic behavior of curvature invariants such as the Kretschmann scalar, and this has been done for the class of polarized Gowdy spacetimes \cite{Chrusciel54} and in more generality \cite{Kichenassamy7}. 
The structure of the horizon and extensions in the polarized Gowdy class can
be very complicated \cite{Chrusciel49,Chrusciel53}.  
In the cases $r$ = 1 ($U(1)$ symmetry) and $r$ = 0 (no symmetry), the large data
global existence and cosmic censorship problems are open. However, in the $U(1)$
case there are conjectures on  the
general behavior which are supported by numerical evidence, and there is a small data semi–-global existence result for the expanding direction  \cite{Choquet01,Choquet03}.

\newpage

\paragraph{Differentiability:}

The technical
questions relating to differentiability are
important from a  mathematical point of view 
regarding well-posedness \cite{Rendall2002}. The
differentiability of the allowed initial data for
the Cauchy problem for a system of partial differential equations 
and the differentiability
properties of the corresponding solutions
are related and determined by the equations themselves. 
For example, in the context of the
Einstein constraints there is a correspondence between the regularity of
the free data and the full data.

There are reasons for considering regularity conditions weaker  than the
natural $C^{\infty}$ condition. One motivation is that
physical matter fields are not necessarily $C^{\infty}$ (so that the theorems need not apply). Another
motivation for considering low regularity solutions is connected to the possibility of
extending (continuing) a local existence result to a global one. 
It is also worth noting that there are 
examples which indicate that generically Cauchy horizons may be non–differentiable
\cite{Chrusciel52}.

There is continued interest in finding a theory for the evolution and constraint
equations for metrics with low differentiability (e.g., to prove the theorems under milder differentiability assumptions such as,
for instance, that the
metric is of regularity $C^{1,1}$  \cite{Kunzinger}
in which the first derivatives of the metric are locally Lipschitz continuous functions,
which is a more natural differentiability class than $C^2$
in a number of physically reasonable situations).
In the existence and uniqueness theorems, the assumptions
on the initial data for the vacuum Einstein equations can be weakened so
that initial data belong to a local Sobolev space.
In spacetime harmonic
coordinates, in which the Einstein vacuum equations form a quasi–linear hyperbolic system,
standard results show that the Cauchy problem is well posed in
an appropriate Sobolev space \cite{HughesKato}, with
improvements on the necessary regularity recently given in \cite{KR,Tataru}.

\paragraph{Singularity theorems:}

The famous singularity theorems are perhaps one of the greatest theoretical
accomplishments in GR and in mathematical
physics more generally \cite{SenovillaGarfinkle}.
Penrose's theorem \cite{Penrose65} was the
first modern singularity theorem, in which the concepts of geodesic incompleteness
(i.e.,  the existence of 
geodesic curves
which cannot be extended in a regular manner within the spacetime and do not
take all possible values of their canonical parameter)
to characterize singularities, Cauchy
hypersurfaces and global hyperbolicity, and closed trapped surfaces \cite{Penrose1979},
were introduced, and has led to many new developments in mathematical GR.
Hawking realized that closed trapped
surfaces will also  be present in any expanding Universe  in its past, which would then inevitability 
lead to an initial singularity 
under reasonable conditions   within GR
\cite{Hawking1966}. 
This subsequently led to the singularity theorem by Hawking and
Penrose \cite{PenroseHawking}, which states that
{\em{if a convergence and a generic condition
holds for causal vectors, and there are no closed timelike curves and there exists at least one
of the following:
a closed achronal imbedded hypersurface,
a closed trapped surface,
a point with re-converging light cone,
then the spacetime has incomplete causal geodesics.}}
It has been argued that
due to the discovery of the cosmic background radiation the singularity theorems
give strong evidence that a singularity actually occurred in our past \cite{HawkingEllis}.

The singularity theorems of Hawking and Penrose proved the inevitability of spacetime
singularities under rather general conditions
\cite{Penrose65,PenroseHawking}. But the singularity theorems say little about the nature of generic
singularities. It should also be pointed out that there are generic spacetimes without singularities
\cite{Senovilla2012}.
For example, the proof of the Penrose singularity theorem does not
guarantee that a trapped surface can arise in evolution. 
Christodoulou  \cite{Christodoulou2009} proved for
vacuum spacetimes a trapped surface can indeed form dynamically from regular
initial data free of trapped surfaces. This result was generalized in
\cite{Klainerman2014} (for more recent work see
\cite{Klainerman2012}). A 
sequence of marginally outer trapped surfaces with areas going to zero which form 
an apparent horizon 
within a region up to the ``center" of gravitational collapse for the 1+3  dimensional Einstein vacuum
equations 
were constructed in  \cite{AnLuk}. 
Marginally outer trapped surfaces also play an important role in proving the positive mass theorem and
the Penrose inequality \cite{AnderssonMetzger} (see below).

There are a number of
open questions, which include proving more general
singularity theorems with weaker  energy conditions and differentiability conditions, and 
determining the relationship
between geodesic incompleteness and curvature 
(e.g., is there always a divergence of
a curvature invariant) \cite{Senovilla2012}.
There are also a number of related open problems in 
cosmology. Generic spacelike singularities are traditionally referred to as being cosmological singularities (but it is not clear that this is necessarily their natural physical interpretation \cite{SenovillaGarfinkle},
since oscillatory singularities might also be related to
the spacelike part of generic black hole singularities \cite{Senovilla2012}; for example,
there is evidence that the mass inflationary instability at the inner horizon of an accreting, 
rotating black hole is generically followed by oscillatory collapse to a spacelike singularity 
\cite{Hamilton}). There is also the question of
singularity resolution in GR by quantum effects and the possibility of singularity theorems in higher dimensions. We shall return to these questions later.

Perhaps the most important open problem within GR is cosmic censorship.

\newpage

\subsection{Cosmic censorship hypothesis}

The Hawking-Penrose theorem  \cite{Penrose65,PenroseHawking} implies that singularities
exist.  
But although the well known Schwarzschild
spacetime  contains a singularity, it is inside the black hole event
horizon and is consequently not visible to outside observers. This leads to the question
of whether gravitational collapse of  realistic matter produces singularities that are similar to the singularity
of Schwarzschild \cite{Penrose1979}, in that they are hidden inside black hole event horizons (weak cosmic censorship)
and are non-timelike (strong cosmic censorship).

Penrose proposed  \cite{Penrose1969} the cosmic censorship hypothesis,
which roughly states that for Einstein's equations 
coupled to ``physical"  matter,
no ``naked singularity" will develop ``generically" from nonsingular
``realistic" initial conditions (Cauchy data). 
A naked singularity is
essentially one with the property that light rays from points arbitrarily
near it can escape to infinity. These singularities are much more disturbing from a physical
point of view, and the question cosmic censorship effectively asks is whether the future can be theoretically
predicted \cite{HawkingEllis}.
It cannot be conjectured that naked singularities never occur, since
there are known examples. However, these examples are of
high symmetry and it is conceivable that naked singularities tend to become
clothed by horizons under most small perturbations.  Indeed, recent results \cite{IsenbergMoncrief1982} tend to
support the notion that naked singularities imply symmetry.

Naked singularities are known to exist in Taub-NUT spacetime \cite{MisnerTaub,Chrusciel53} and
simply by removing regions from Minkowski spacetime.
It is also known that the
equations of a pressureless fluid or ``dust" will lead to spurious ``shell
crossing" naked singularities. 
In particular, a central locally naked singularity forms in spherical dust Tolman-Bondi-de Sitter collapse \cite{Goncalves} from a non-zero-measure set of regular initial data, 
at which the Weyl and 
Ricci curvature scalars diverge.
The most comprehensive results known on global inhomogeneous solutions of the Einstein
equations are for solutions of the spherically symmetric 
Einstein equations coupled to a massless scalar field with asymptotically flat
initial data, where Christodoulou
has proved  that naked singularities can develop from
regular initial data \cite{Christodoulou94} and that this phenomenon is unstable with respect to
perturbations of the data \cite{Christodoulou99b}.

Consequently, we seek to formulate cosmic censorship as a precise mathematical
conjecture and then find a proof or a counterexample. 
Theorems on maximal Cauchy developments are within the 
global theory of partial differential equations and are generally very
difficult to prove  \cite{Ringstrom2010}.
There can be no timelike singularities in a globally hyperbolic spacetime. Thus,
a method for formulating (strong) cosmic censorship is as a statement that (under suitable
conditions) spacetime must be globally hyperbolic. However, an initial data set has a maximal
Cauchy development, which is a globally hyperbolic spacetime, but that maximal
Cauchy development may not be the complete spacetime.

There are two other particular problems that must be faced.
First, a naked singularity is very difficult to accurately define mathematically. Since the Einstein equations are essentially hyperbolic, the notion of extending a solution to points which can ``see" the singularity is problematic, and so we have to seek an alternative definition
of a naked singularity that is more stable and can be
mathematically formulated. 
The second problem is genericity. It is known that there are special examples of solutions in GR which,
for all reasonable definitions, contain a 
naked singularity where the maximal
development is extendible.
So it is impossible to prove a general statement that 
says a naked singularity cannot exist. That is, without some sort of ``generic condition'', this version of
cosmic censorship would fail. We are, of course,
ultimately interested in the real process of
gravitational collapse, but care must be taken not to formulate a 
conjecture that will be vulnerable to what a physicist might claim appears to be
an artificial counterexample.  Therefore, the aim is to
refine the conditions of the conjecture to rule out non--physical
counterexamples, but not to the extent of making cosmic censorship irrefutable \cite{Rendall2002}.

There are actually two different cosmic censorship hypotheses, which are only minimally 
related to each other. The weak cosmic censorship hypothesis states that:
{\em{For generic initial data to the evolution problem in GR, there cannot be naked singularities.}}
This is such an open problem that the correct formulation of the statement is not even known  \cite{Christod99}.
For an extensive treatment (including a somewhat
precise version) of the weak cosmic censorship conjecture see \cite{Wald1998}.

\vspace{3pt}

\noindent
{\bf{Problem P1:}} Prove the weak cosmic censorship conjecture.

\vspace{3pt}

In the case of asymptotically flat spacetimes (describing isolated systems
in GR), the
work of Christodoulou establishes weak cosmic censorship in the class of spherically
symmetric Einstein–-scalar field spacetimes \cite{Christod}, and also gives examples of initial
data such that the Cauchy development has a naked singularity \cite{Christodoulou99b}.

The second hypothesis is strong cosmic censorship, which
states that:
{\em {A generic solution to the Einstein's equation cannot be continued beyond the Cauchy horizon.}}
For earlier
surveys on the strong cosmic censorship conjecture, see \cite{Isenberg92,Chrusciel49}.
It is of interest to prove weak and strong cosmic censorship even for
vacuum solutions of Einstein's equations (i.e., those with no matter) or,
more generally, within special classes of spacetimes.
\vspace{3pt}

\noindent
{\bf{Problem P2:}} Let M be a 3 dimensional compact manifold. Prove that the maximal
vacuum Cauchy development for generic vacuum data sets  is equal to the maximal vacuum
extension of M.

\vspace{3pt}

\newpage

An alternative strategy is to search 
for a counterexample  to cosmic censorship. If a wide class of possible counterexamples can be shown to fail,
then this might even be seen as evidence for the likely validity of the conjecture. 
A possible counterexample for
weak cosmic censorship might arise  
from a process in which a
black hole  turns into a naked singularity.  For example, the Kerr metric with mass M and angular momentum
J represents a black hole if $J\leq M^2$ and a naked singularity if $J > M^2$. Therefore, 
a naked singularity might possibly be produced by overspining a black hole. 
Since  spinning black holes repel
the particles whose angular momentum would increase their spin, such a ``spin-spin
repulsion" unfortunately prevents the overspinning of a black hole \cite{SenovillaGarfinkle}.

A plausible candidate
for a vacuum counterexample to cosmic censorship (with a negative cosmological constant) has recently been proposed
based on the superradiant instability of Kerr-–AdS black holes \cite{Dias2015}.
Another plausible counterexample 
(based on a holographic
model of an electrically charged localised defect)
in four-dimensional
Einstein-–Maxwell theory with asymptotically anti--de Sitter boundary
conditions was presented in
\cite{Horowitz2016}; smooth initial data was shown to evolve to a region of arbitrarily large curvature
in a finite time that is visible to distant observers. 
Unlike the spherical collapse `‘counterexamples’'
which are finely tuned, this candidate is generic \cite{Bizon}.

%\newpage 

Finally, we note that by 
considering only globally
hyperbolic spacetimes,
solutions with
gross causality violations are excluded, while some singular behaviour is still possible. 
But there are exact solutions with closed
timelike curves known (e.g., the Godel and NUT spacetimes). 
The existence of such causality violation gives rise to ``existential
problems of an imponderable nature'' \cite{HawkingEllis}.
Stephen Hawking has
suggested the ``chronology protection conjecture" 
that asserts that the closed timelike curves which arise in some solutions to the equations of GR (and
which imply the possibility of backwards time travel) will be ruled out by a future theory of quantum gravity.

\newpage

\subsection{Penrose inequality}

The mathematical ideas behind the proofs of the singularity
theorems have been applied to several important results in GR, such as the positive
mass theorem in its original form \cite{SchoenYau} which has, in turn,  led 
to research on the rigidity of asymptotically flat manifolds with non-negative 
scalar curvature. 

In particular, Penrose has shown \cite{Penrose1973} that if a certain inequality
involving the area of a marginally (outer) future-trapped surface
(the apparent horizon) and the (ADM) mass of the initial hypersurface containing this horizon were
violated, then the spacetime that results
from evolving the initial data contains a naked singularity. Therefore, initial data
violating this so-called {\em{Penrose inequality}} would constitute a counterexample to weak cosmic
censorship, while a proof of this inequality would provide evidence in favor of
weak cosmic censorship. In fact, such a proof
would possibly lead to an approach for attacking the cosmic censorship conjecture
using methods in partial differential equations \cite{Moncrief81E} (this is
discussed further in  \cite{LARS99}).

\vspace{3pt}

\noindent
{\bf{Problem P3:}} Find a proof of the Penrose inequality or present a counterexample in
the general case.
\vspace{3pt}

The Riemannian version of the Penrose inequality was recently proved \cite{Huisken}.
The proof in the Lorentzian case is not known. Even in spherical symmetry only a weaker version (using the
energy rather than the mass) is known to hold.
 Proofs have been given
under various restrictive assumptions, such as the existence of certain
foliations (e.g., the constant mean curvature time gauge \cite{Moncrief96}), and global conditions on the spacetime
(see \cite{Bray}).

The  Penrose inequality is one of a large class of mass inequalities for spacetime manifolds
\cite{DainS}; for example, an analogous inequality is based on the Penrose
quasi-local mass \cite{Tod}. It is also of interest to find a generalisation
of the Penrose inequality to initial data sets which are not time-symmetric.
There also exist stronger versions of the Penrose inequality involving angular
momentum, electric charge, and/or the cosmological constant \cite{SenovillaGarfinkle},
most of which lead to  open questions. 
There are further refinements of the conjectures, such as the so-called Gibbons-Penrose inequality, which 
gives some improved lower bounds when there are multiple black holes \cite{Gibbons}. 
Another inequality is Thorne's hoop conjecture
\cite{Thorne1972}, which exploits
the physical idea that since black holes are extremely localized objects, 
their energy/matter content must be severely compacted in all spatial directions.
Despite the difficulty in 
making this idea  precise, the hoop
conjecture has proven successful \cite{Wald1998}.
A possible mathematically viable
reformulation of the conjecture has been presented in \cite{Senovilla2008}.

%\newpage

%----------

\newpage

\subsection{Yang--Mills equations and GR}

Many important mathematical
questions, including 
stability theorems and the proof of existence of
Yang-Mills (YM) fields by methods of partial differential equations, remain unsolved.
We have already discussed problem M2 on the existence of solutions of YM  earlier, and there
were 2 well known problems (Y117, Y118) presented 
in \cite{Yau1982}, the first of which is the question of whether every SU(2) Yang--Mills field is self-dual or anti--self--dual.
A key contemporary
question is to obtain solutions of the YM equations on a Riemannian (or Lorentzian)
manifold. Recently it has been shown numerically that the static, spherically
symmetric Einstein-Yang-Mills (EYM) equations have non-singular, asymptotically flat
solutions 
\cite{BartnikMckinnon}.
Six interesting questions for EYM solutions
were presented in \cite{Bartnik} (RB17, see above).

A central feature of YM theory  is  the invariance
of the physics under an infinite-dimensional
group, in which bundles, connections and curvature
play a fundamental role. It is consequently 
a subject of interest not only to physicists but also,
particularly  after the work of Atiyah,
Hitchin, and Singer \cite{AHDM,AHS1,AtiyahSinger}, to mathematicians (as discussed earlier and in
\cite{Yau1982}).
The YM field equations depend on
how a section of the Lie algebra valued bundle is choosen. The choice of
such a section is called the choice of a gauge. In a suitable gauge, the
YM equations become a quasilinear elliptic system whose highest
order term is linear.
Physicists are mostly interested in YM fields over $ R^4$ or $S^4$.

It is known that in four dimensions there exist global smooth solutions of the
YM equations corresponding to rather general initial data. 
Global existence in Minkowski space, assuming initial data of
sufficiently high differentiability, 
was first proven in \cite{Eardley} and a new proof
of a local existence theorem for data of finite energy (and since energy
is conserved this immediately proves global existence) was given in  \cite{Klainerman95}.
A global existence proof  on 1+3 dimensional, globally hyperbolic
spacetimes was given in \cite{Chrusciel56} (see also \cite{Kichenassamy7}).
The proof of the global existence to the future for
hyperboloidal initial data  close to the standard data on a hyperboloid in Minkowski
space by Friedrich was
later generalized to Maxwell and YM matter in \cite{Friedrich91}. 
However, although asymptotically flat (with regular interior) spherically symmetric and localised 
(``particle-like") solutions of the coupled EYM equations with gauge group SU(2) have 
been known for many years, their properties are still not well understood \cite{Olinyk}.

In dimensions greater than five it is known that there exist solutions which develop
singularities in a finite time.
Numerical evidence indicates that this type of blow-up is
stable (i.e., it  occurs for an open set of initial data) and
that there is a critical self-similar solution separating this kind of blow-up
from dispersion. There is as yet no rigorous proof of blow-up in five dimensions.
In six dimensions
singularities form, but apparently differently from those in five
dimensions \cite{LARS99}.

\newpage

The effects found in YM theory are captured in two dimensions less
by {\em{wave maps}} with values on spheres, where it is easier
to prove theorems.
The existence of a solution
having the properties expected of the critical solution associated with singularity formation for wave maps in four dimensions has been proven in \cite{Bizonwavemaps}.
An important open question is the global
existence problem for the classical  wave map equation (i.e., the nonlinear $\sigma$-–model, hyperbolic
harmonic map equation).
The wave map equation has
small data global existence for spatial dimension $n \geq 2$. 
But global existence for “large data” is known only for symmetric
solutions and, in particular, the global existence problem for the wave map
equation is open for the case $n = 2$. For the case $n = 1$, global existence
can be proved using energy estimates  \cite{LARS99,AnderssonGudapatiSzeftel}.
The U(1) symmetric vacuum 1+3 case in which 
the Einstein equations reduce to
1+2 gravity coupled to wave map matter in the presence of a hypersurface
orthogonal space–like Killing field, 
is of intermediate difficulty between the full 1+3 Einstein equations
and the highly symmetric Gowdy equations \cite{SterbenzTataru}.

It is also of interest to consider other forms of matter such as,
for example,  self-gravitating
collisionless matter models (see the reviews \cite{Rendall2002,LARS99,Andreasson} and references within).
There are theorems on the global existence and uniqueness of
smooth solutions of the
Vlasov-Poisson  and the classical Boltzmann equations 
in Newtonian theory. Many
analogues of these results 
have been proven in GR, including the
global existence of weak solutions, the
convergence to equilibrium for classical solutions starting close to equilibrium,
basic existence theorems for spherically symmetric static solutions, 
plane and  hyperbolic symmetric spacetimes and
a subset of general Gowdy spacetimes, and studies of spherically symmetric collapse. 
Collisionless matter models
are known to admit a global singularity-free evolution, and in
many cases can
also lead to isotropization at late times.
Analytical techniques have not been applied in the general case, although 
numerical methods have been used to gain some insights \cite{Rendall2002,LARS99}.

\vspace{3pt}

\noindent
{\bf{Problem P4}}: Prove the global existence of classical  spatially inhomogeneous solutions for
small initial data in collisionless matter models. Prove an existence and
uniqueness theorem  for general spatially homogeneous
(such as Bianchi type IX) solutions of the Einstein-Vlasov equations and investigate the
large initial data case.

\vspace{3pt}

%%%%%%%%%%%%%%

\newpage

\subsection{Uniqueness and stability}

Mathematically, any proof of stability requires deriving the asymptotic behaviour of solutions to the Einstein equations in GR, a highly nonlinear system of partial differential equations, which is 
notoriously very difficult. However, there are  some special cases for
which there exist proofs or which have received particular attention. 

\paragraph{Stability of Minkowski spacetime:}

Minkowski spacetime is globally stable \cite{Christodoulou93}. That is, if we start with a universe that
is already very sparse, it is guaranteed that it will evolve asymptotically to Minkowski spacetime. 
The first result on  the global existence (for small data) and the stability of Minkowski spacetime 
under the field equations of GR was due to Christodoulou and Klainerman  \cite{ChristodoulouKlainerman90,Christodoulou93}. They proved that if initial data for the vacuum Einstein equations
are prescribed which are asymptotically flat and sufficiently close to those
induced by Minkowski spacetime on a hyperplane, then the maximal Cauchy development of this data is geodesically
complete (and they further  provided details on the asymptotic behaviour
of the solutions). 
Results can also be found for any asymptotically flat spacetime
where the initial matter distribution has compact support, so long as attention
is confined to a suitable neighbourhood of infinity.
There are recent extensions to these results
by various authors (e.g., see  \cite{Bieri}).

\paragraph{Uniqueness of black holes:}
If we conjecture that the final state of a spacetime is either Minkowski space or a black hole, 
we can then ask whether a black hole 
is the only possible stationary  (steady state) solution. The problem of black hole 
uniqueness is  not completely resolved. The study of uniqueness for non-vacuum spacetimes 
is colloquially known as  ``no-hair"
theorems.

In the case where it is assumed that the spacetime has additional symmetry  and is either axially symmetric or rotationally symmetric, 
the uniqueness of black holes is known. 
The  uniqueness of the  4D Schwarzschild and Kerr solutions in GR was discussed in \cite{HawkingEllis}.
The uniqueness theorem for Schwarzschild spacetime was presented in \cite{Bunting,Israel}.
The unique stationary (nonstatic) regular predictable
Ricci flat spacetime subject to certain assumptions is 
the Kerr solution \cite{Kerr}.
The uniqueness theorem for the Kerr spacetime was proven in \cite{Carter,Robinson}.
In the non-vacuum case  the uniqueness of the rotating electrically charged black hole solution
of Kerr-Newman has not yet been generally proven (however, see  \cite{Newman65,Mazur}).

We also know that black holes are unique 
if we assume real analyticity.
If the regularity assumption is relaxed to just infinitely differentiable the result is still expected to be true.
In this case there are only some partial results. For example, if
only small perturbations of a stationary black hole are allowed then there are no other stationary solutions that are approximately a known 
black hole solution without being one, and if  certain special structures on the event horizon
are assumed then other stationary 
exteriors are not possible.

\newpage

\paragraph{Stability of Kerr-Newman black hole:}
If we assume that the known Kerr-Newman family of black holes form the unique stationary state of GR, 
the next problem is to prove that they are actually stable under perturbations. That is, if we start out with initial data very close to that of a Kerr-Newman black hole, does the the evolution ``track" a Kerr-Newman black hole.
Although there has been substantial and exciting progress made in the linearised problem \cite{DHR}, results for the full nonlinear problem are still elusive.

The stability of the Kerr metric was discussed in \cite{Heusler}, and a comprehensive review 
was given in  \cite{Chand}. The aim is to show that
perturbations of the
Kerr (and Schwarzschild \cite{Holzegel1016}) solution decay exponentially and are thus stable.
Unfortunately,
a mathematically rigorous understanding of the stability of
the generic Kerr black hole, as well as a thorough understanding of
its dynamics under arbitrary non-linear perturbations, is still lacking. 
However, current observational
data are compatible with the predictions of GR,
and suggest that the end point of mergers is a Kerr
black hole. Indeed, all numerical results provide evidence
that the Kerr (and Kerr-Newman) black holes are non-linearly stable (at least
within a certain range of the angular momentum)
\cite{Zilho}.

\vspace{3pt}

\noindent
{\bf{Problem P5:}} Prove the stability of the Kerr black hole.
\vspace{3pt}

It is of interest to extend stability results to the case of a non-zero cosmological constant
\cite{Dotti}.
Regarding the  stability of the Kerr-de Sitter
family of black hole solutions, there has been recent results on nonlinear perturbations
in the slowly rotating case \cite{Hintz}.
The case of a negative cosmological constant is much more problematic, because it is not even clear if 
the Kerr-AdS black hole is itself stable (due to superradiance and stable trapping
phenomena \cite{Holzegel}).
We shall discuss the stability of the de--Sitter and anti--de--Sitter spacetimes later.
It is also of interest to study the stability of models with matter, particularly in the cosmological context (also see later). Unfortunately, even generalizations to simple inhomogeneous perfect fluids are problematic
since the formation of shocks (or, in the case of dust, shell-crossings) are anticipated
to occur which form a barrier to the mathematical study
of the evolution of the cosmological
models with known techniques.
Criteria for the development of shocks (or their absence),
based on the techniques of classical hydrodynamics, should be developed further.

\newpage

\subsection{Other problems}

\paragraph{Curvature invariants:}
In \cite{CHP} it was shown that the class of 4D Lorentzian manifolds that cannot
be completely characterized by the scalar polynomial curvature invariants constructed from
the Riemann tensor and its covariant derivatives must be 
of a special ``degenerate Kundt form''. This result, which is also believed to be 
true in higher dimensions \cite{CHP2},
implies that
generally a spacetime is completely characterized by its scalar curvature invariants
(at least locally, in the space of Lorentzian metrics).
The special Kundt
class is defined by those metrics admitting a null vector that is geodesic, expansion-free,
shear-free and twist-free.
We recall that in the Riemannian case a manifold is always locally characterized by its scalar polynomial
invariants.

It is also of interest to study (the `inverse question') of 
when a spacetime can be explicitly constructed from its scalar curvature invariants.
In 4D we can (partially)
characterize the Petrov type of the Weyl tensor in terms of scalar curvature invariants \cite{krameretal}. 
Having determined when a spacetime is
completely characterized by its scalar curvature invariants, it is also of interest to determine
the minimal set of such invariants needed for this characterization.
\vspace{3pt}

\noindent
{\bf{Problem PF1}}:
Determine
when a 4D spacetime can be explicitly constructed from its scalar curvature invariants and
determine the minimal set of such invariants.
\vspace{3pt}

\paragraph{Evolution of the horizon:}  
There is much interest in determining the appropriate definition of the ``boundary of a black hole". 
A closed oriented space-like 2-surface (normally
isomorphic to $S^2$) in a spacetime determines two future
null vector fields, normal to the surface. If the future evolutions of the surface
along these directions are both area-non-increasing, the surface is future
trapped, and if one of the null mean curvatures
is zero, then the surface is called an ``apparent horizon''. 
It is also  important to determine the  evolution of the horizon and, more generally, formulate an appropriate
definition of a dynamical horizon in GR.
We note that much work
on the evolution of apparent horizon (such as black hole evaporation) is based on a linear analysis, which to first order assumes 
that the horizons do not move. The true nonlinear versions of the evolution is not yet well understood.
The problem of 
identifying and locating horizons using scalar curvature invariants has recently been studied \cite{ColeyMcNutt}.

\vspace{3pt}

\noindent
{\bf{Problem P6:}} Formulate an appropriate
definition of a dynamical horizon.

\vspace{3pt}

\newpage

\paragraph{Geodesic hypothesis:} One of the postulates of GR is that point particles with negligible mass will travel 
along geodesics of the spacetime. 
\vspace{3pt}

\noindent
{\bf{Problem P7}}: Prove that test particles  move on spacetime
geodesics.
\vspace{3pt}

\noindent
This famous problem (RB23) was first considered by Einstein in the 
1920s and is still not completely resolved (there has been an extensive investigation
using asymptotic expansions --- see the discussion in \cite{JEhlers}).
The main problem is how to make the process of ``taking the negligible mass limit" 
rigorous. And for a physical object in GR, when it moves, its motion will cause ``ripples" in the spacetime caused by 
gravitational backreaction of its own presence. 
In addition, while the three body problem is difficult in classical mechanics, even the two 
body problem in full generality is still unresolved in GR.

\paragraph{Newtonian limit:}
It is difficult to give a precise mathematical formulation of
the statement that Newtonian gravitational
theory is the limit of GR as the speed of light tends to
infinity. 
Ehlers gave a definition of the Newtonian limit of GR which encodes
those properties which are physically desirable \cite{Ehlers}. However, even when a suitable 
definition has been given, the question still remains as to whether the definition is compatible
with GR in the sense that there are general families
of solutions of the Einstein equations which have a Newtonian limit with
the chosen definition. Asking whether there are such families which are suitably
differentiable  is related to the issue of giving
a mathematical justification of the so--called post-Newtonian approximation. 
See problems RB20 and RB21.

\newpage

\section{Theoretical physics problems in the Quantum realm (Ph1-Ph4)}

There are a number of fundamental questions in the quantum realm, culminating in the ultimate question of whether
there is a single theory (or even, more precisely, one single equation)
that would unify all of nature within a so-called ``theory of everything".
In particular, is this theory string theory?
And would this theory then give an explanation of the fundamental gauge group in the
grand unification theory of the three non--gravitational forces \cite{Arkani}
and also
explain the values of all fundamental physical constants (and 
whether they vary over time)?  
In addition, are there fundamental
particles that have not yet been observed and, if so, what are their properties?  
Let us consider the following particular problems.

\paragraph{The foundational problems of quantum mechanics (Ph2):}
These problems concern the fundamental understanding of quantum physics 
and especially the important role that measurement
and observation play in the description of physical reality.
There are currently many
interpretations of quantum physics, including the classic Copenhagen interpretation, Everett's controversial ``many
worlds'' interpretation, and even more controversial ones such as the ``participatory anthropic principle''.

In particular, how does the quantum interpretation of
reality, which includes the superposition of states
and wavefunction collapse or quantum decoherence, give rise to what we perceive?  
What are the actual causes of the collapse of the quantum
wavefunction?  
Are there non-local phenomena in quantum physics 
and, if they do in fact exist, are they limited to the
entanglement revealed in the violations of the Bell inequalities and   
can they be observed?  
What does the existence or
absence of non-local phenomena imply about the fundamental structure of spacetime and how is this related to quantum
entanglement?  
Most modern physicists who work within quantum field theory perhaps no longer consider
questions of the proper interpretation of the fundamental nature of quantum physics to be
of prime importance.  Indeed, many may believe that the principle of decoherence is essentially an appropriate explanation; 
for example,  interaction with the environment causes the
quantum collapse.

%%\newpage

However, dynamical models have been proposed to explain the collapse of the wave-function and
perhaps provide a possible solution
to the quantum measurement problem, by proposing that the Schrodinger equation is an approximation to a stochastic nonlinear dynamics (with the stochastic
nonlinear aspect becoming increasingly more important when progressing from microscopic systems
to macroscopic ones)  \cite{TejinderSingh}. 
In addition, as in most other physical systems, evolution in time is central to the understanding of
quantum systems. The time that is used to define evolution in quantum theory is clearly
part of the classical spacetime manifold. However, this perhaps suggests that the
present formulation of quantum theory is incomplete and that there ought to
exist a reformulation of quantum theory which does not refer to classical time.

%%\newpage

\paragraph{The unification of particles and forces (Ph3) and the tuning problem (Ph4):}
The standard model of particle physics involves eighteen different fundamental particles.
It is often believed that a theory of nature should have a more fundamental method of unifying these
particles.  For example, string theory, which is perhaps the most well-defined
approach, predicts that all particles are different vibrational modes of fundamental filaments of energy or
strings. It is, of course, of great importance to
determine whether or not the various particles and forces can be unified within a theory that explains them all as manifestations of a single, fundamental entity.

In the standard model of particle physics the parameters representing the eighteen particles
predicted by the theory are  required
to be determined (i.e., measured by observations) in order for theoretical predictions to be made.
However, some physicists argue that fundamental physical principles of a unified field theory should set these
parameters, independent of measurement. In particular, there is the question of whether the form of the universe
is inherently set by its properties (in the sense that the properties would not occur if the
form is different).
In the multiverse paradigm there is not just a single universe, but
there are a wide range of fundamental theories (or different variants of the same theory, based on different physical parameters) and our universe is just one of the possible universes that could be created.
In this paradigm the question then becomes why our 
particular universe has properties that appear to be so finely tuned to allow for the
existence of life. This has led some scientists to turn to the
anthropic principle for explaining this fine-tuning problem:
this asserts that our universe must have the properties it does because if it had
different properties we wouldn't be here to be able to beg the question. 

Finally, 
the recent observation of a Higgs boson appears to complete
the standard model, but with the addition of new physics in order to protect
the particle mass from quantum corrections that would increase it by perhaps fourteen
orders of magnitude or more. It is widely thought that the most plausible resolution of this
hierarchy (or naturalness) problem is
supersymmetry. However, the simplest supersymmetric models have
not worked, and no convincing mechanism has yet been found to either
break supersymmetry or to determine the many supersymmetric parameters (AL3.3).

%%\newpage

\paragraph{The problem of quantum gravity (Ph1):}
There are four fundamental forces of physics. However, the standard model of particle physics includes only 
the three forces of 
electromagnetism and the strong and weak nuclear forces.  
An attempt to formulate a theory which unites all four forces,
including gravity,  into a single unified field theory is a primary goal of
theoretical physics.
The theory that includes both GR and
the standard model of particle physics is referred to as quantum gravity (QG).  
Unfortunately, at present these two theories describe different scales of nature and any
attempts to explore the overlaping scale has yielded incomprehensible results, such as the force of
gravity (or curvature of spacetime) becoming infinite.  If quantum mechanics and GR can be realized 
as a fully consistent theory of QG a number of natural questions arise, which include
the following:  Is spacetime fundamentally continuous or discrete?  Would such a theory include a
force mediated by a hypothetical graviton, or would it be a product of a discrete structure of spacetime itself 
(such as, for example, in loop
quantum gravity)?  Are there possible differences from the predictions of GR at very small or very large
scales (or in any other extreme circumstances) that result from a theory of QG?

Although QG effects modify GR, leading to new
gravitational physics, it appears that these modifications do not significantly affect the macroscopic
behavior of stellar systems and black holes. For example, 
a black hole that
evaporates through the emission of Hawking radiation \cite{Hawking},
perhaps the most dramatic consequence of uniting
GR and quantum mechanics, does not differ significantly from a
classical black hole over astrophysical timescales \cite{Barcelo}.

%%%%%%%%%%%%%%%%%%%%%%%%%%%%

Many of the above problems are in theoretical physics and generally are not problems in mathematical
physics. However, a lot of mathematics is utilized in
{\em{string theory:}} for example,
Yau proved the Calabi conjecture, which allowed physicists to show, using Calabi–-Yau compactification, that string theory is a viable candidate for a unified theory of nature. In addition, 
there are a number of related  fundamental questions in 
Yang-Mills theory, which we have discussed earlier.
We shall next discuss some specific problems that are definitely within mathematical physics.

\newpage

\subsection{Instability of Anti-de Sitter spacetime}

Anti-de Sitter (AdS) spacetime in any dimension is the unique maximally symmetric Lorentzian manifold
with constant negative scalar curvature. 
AdS spacetimes are of interest in theories of QG formulated in terms of string theory (in which elementary particles are modelled not as zero-dimensional points but as one-dimensional objects called strings) or its modern extension, M-theory. Indeed, AdS spacetimes have
come to play a central role in theoretical physics, primarily due to the AdS/CFT
correspondence (or Maldacena gauge/gravity duality)
which is the conjectured equivalence between string theory on
an asymptotically AdS spacetime and a conformally invariant quantum field
theory (CFT) living on the boundary of this spacetime \cite{mald,Klebanov}.
CFT are quantum field theories, including theories similar to Yang–-Mills theories, that describe elementary particles.

The AdS/CFT correspondence suggests that it is possible to describe a force in quantum mechanics (like electromagnetism, the weak force or the strong force) in a certain number of dimensions with a string theory where the strings exist in an AdS spacetime with one additional dimension.
The duality represents a major advance in our understanding of string theory and QG
since it provides a non-perturbative formulation of string theory with certain boundary conditions.
The usefulness of this strong-weak duality results from the fact that 
strongly coupled quantum field theories can be studied by investigating the corresponding
weakly interacting gravitational theory which is mathematically more tractable. This has been used to study many aspects of nuclear and condensed matter physics 
(such as, for example, the modelling of
non-equilibrium processes such as heavy ion collisions)
by translating those problems  into more mathematically tractable problems in string theory.
That is, the AdS/CFT dictionary is used to translate 
the strongly coupled
CFT to the string dual, which effectively reduces to classical AdS gravity, 
and the results are then utilized to produce useful information on the 
physics of the CFT. Unfortunately, the problem with this holographic approach  is that
the gravity side  in the non-stationary regime is not well understood.

\newpage
The AdS/CFT correspondence
provides strong motivation for studying the dynamics of asymptotically AdS
spacetimes. But, of course, this
is an interesting problem in classical GR in its own right.
AdS spacetime is different to 
Minkowski and de-Sitter spacetimes \cite{HolzegelShao}, which were proven to be nonlinearly stable a long time ago
\cite{Christodoulou93,Friedrich}.
It has recently been conjectured that the AdS spacetime is unstable
under arbitrarily small perturbations  \cite{PiotrBizon}. This is related to some interesting more general
mathematical problems.

The question of the global
nonlinear stability of AdS was given a huge boost by the seminal work of Bizon
and Rostworowski \cite{Bizon} following a conjectured instability by Dafermos and Holzegel
\cite{Dafermos}. While it would of course be desirable to study the nonlinear stability of AdS with no
symmetry restrictions, this problem currently is analytically and numerically intractable.
In \cite{Bizon} the analysis was restricted to spherical
symmetry within pure Einstein
gravity with a massless scalar field.
Numerical results  suggested that AdS is nonlinearly unstable to a weakly
turbulent mechanism that forms an arbitrarily small black hole, whose mass is controlled
by the energy of the initial data. While this
nonlinear instability seems to occur for generic perturbations, there are perturbations that
do not necessarily generate an instability  (see \cite{DiasSantos} and references within) which, in turn, 
appears to lead to the
existence of islands of stability \cite{DiasGreen}. Using
standard perturbation theory to third order in the amplitude of the linear seed, it was shown \cite{Bizon} that this
leads to secular growth and
nonlinearities occur that can
create resonances. 
The heuristic explanation for the mechanism which triggers the turbulent behaviour 
is thus the generation of secular terms by resonant four-wave interactions; it is this weak turbulence that is a driving mechanism of the instability.

There are modifications of standard perturbation theory that can capture the dynamics up
to certain time scales, such as the resonant approximation \cite{Bizon2015}, but rely on the 
spherical symmetry assumption.
It is not known if any solution of the Einstein equations with a fully resonant
spectrum necessarily possesses a nonlinear instability, but it is clear it is a necessary condition
for the existence of the weakly turbulent instability. 
It is an interesting open question as to whether the nondispersive character of the linearized spectrum
is essential for the
turbulent instability and how generic is the turbulent
instability. 
In order to study this beyond spherical
symmetry 
third order perturbation theory calculations for a variety of different seeds have been performed  \cite{DiasSantos,DiasGreen}, and it was found
that the gravitational case is more richer than the spherically symmetry case analysed
in \cite{Bizon}. The prime question is consequently to determine the endpoint
of instability of arbitrary dimensional AdS spacetimes for non-spherical perturbations  \cite{PiotrBizon}.
Note that recently nonlinear instability was
proved for the spherically symmetric Einstein-massless Vlasov system
\cite{Moschidis}.

\vspace{3pt}

\noindent
{\bf{Problem P8:}} Determine 
whether the conjectured nonlinear instability in Anti-de Sitter spacetime, which leads to a weakly 
turbulent
mechanism that develops a cascade towards high frequencies leading to black
hole formation, behaves differently in more general models than
spherically symmetric scalar field collapse.

\vspace{3pt}

\newpage

%%%EDIT FURTHER--OMIT?

Since this is a particularly topical problem, let me discuss it in a little more detail \cite{Anderson}.
In the case of AdS, the question of stability must be supplemented by a
choice of boundary conditions at infinity and, a priori, any results may depend
on this choice.
A local well-posedness result to the initial boundary
value problem for
a large class of AdS boundary conditions was proven by Friedrich for the vacuum
Einstein equations with negative cosmological constant in 4D
\cite{Friedrich}, allowing local stability to be studied mathematically. 
In the case of  reflective boundary conditions, for which
there is no flux of energy across the conformal boundary, the
asymptotic stability of AdS is not possible because the (conformal) boundary acts
like a mirror at which perturbations propagating outwards bounce off and return
to the bulk. This leads to very complex nonlinear wave interactions
in the bulk, which is extremely difficult to study even in the case of small perturbations.
Consequently, it is hardly surprising that the question of the stability of AdS spacetime
remains open.

For reflective boundary conditions, the
problem of the linear stability of AdS
reduces to a much simpler spectral problem for a certain master linear
operator whose coefficients depend on the character (i.e., scalar, electromagnetic
or gravitational) of the perturbations \cite{Breitenlohner}.
The problem of the nonlinear  stability of $n+1$ dimensional AdS spacetime
in full generality is currently beyond the theory of partial differential equations. 
Thus it is natural to consider more tractable special
cases.  In particular, for  spherically symmetric perturbations of a self-gravitating
minimally coupled massless scalar field,
the system of Einstein-scalar field equations with
appropriate boundary conditions and compatible smooth initial data 
constitutes a locally well-posed initial-boundary value problem in asymptotically AdS spacetimes.
Perturbative and numerical studies
of the global behavior of small data solutions to this problem 
give evidence (first for $n = 3$ and later generalized
to $n \geq 3$ \cite{Jalmuzna}) for the conjecture that (within the model)
the $n+1$ dimensional AdS spacetime  is unstable to the
formation of a black hole for a large class of arbitrarily small perturbations \cite{Bizon}.

\newpage

%%%%%%%%%%%%%%%%%%%%%%%%%%%%

\subsection{Higher dimensions} 

Extra dimensions (beyond the
familiar four of ordinary spacetime) are employed in string theory \cite{string}. In addition,
spacetime manifolds of higher dimensions are considered in some cosmological scenarios. 
If there are, in fact, higher dimensions, then deep questions on  the structure
of the internal space for our Universe arises (AL5.1).
If nature has more than four spacetime dimensions, what are their size, what is the topology of Universe,
and why are there 3  apparent spatial dimensions?
And can we experimentally observe
evidence of higher spatial dimensions?

The study of black holes in GR, and the differences
between black holes in 4D and higher dimensions, is
currently of great interest.  At the classical
level, gravity in higher dimensions exhibits a much richer structure than in 4D; 
for example, one of most remarkable features of 4D GR is the uniqueness of the Kerr 
black hole.  In contrast, there exist a number of
different asymptotically flat, higher-dimensional vacuum black hole
solutions \cite{EmparanReall}.  
The uniqueness and stability of higher dimensional black holes is of  paramount interest.

\vspace{3pt}

\noindent
\paragraph{Problem P9:}
Determine the uniqueness of black holes in higher dimensions.

\noindent
\paragraph{Problem P10:}
Determine the stability of  higher dimensional black holes.
\vspace{3pt}

A number of sub--problems,
including adapting the hypotheses of analyticity, non-degeneracy, and connectedness in the
black-hole uniqueness theorems, and
classifying all vacuum near-horizon geometries with compact cross-sections, have been proposed.

Differential geometry and geometric results have recently been developed
in higher dimensions \cite{CHP,higherdimensions}. In particular, 
even though the singularity theorems were originally proven in 4D, results 
in which the closed trapped surface is a co-dimension
two trapped submanifold hold in arbitrary dimensional spacetimes. 
The concept of being trapped can also
be associated with submanifolds of any co-dimension,  
so long as
an appropriate curvature condition is assumed to ensure the existence of focal points to the submanifold.
Closed trapped surfaces
in co-dimension 3 and in arbitrary dimension were discussed in \cite{Galloway}. 
The positive mass theorem has recently been proven in all dimensions \cite{SchoenYau2017}.

\newpage

The question of stability in higher dimensions is more problematic.
On one hand, radiative decay of solutions is stronger in higher dimensions and would enhance stability. On the other
hand, 
there are more degrees of freedom which will generally increase the possibilities of instability. There is numerical evidence to suggest that 
certain types of higher dimensional black holes are in fact unstable \cite{EmparanReall}.

However, the problem of
cosmic censorship in higher dimensions is not well posed and very difficult. 
Indeed,
in higher dimensions
there is strong numerical evidence that cosmic censorship fails \cite{LehnerPretorius2010}, and
higher dimensional black holes can be unstable under
gravitational perturbations. This was first shown by Gregory
and Laflamme for black strings and black p-branes \cite{GregoryLaflamme} (in 4D this
instability does not exist).
Convincing numerical 
evidence that unstable black strings pinch off in finite
asymptotic time, thus resulting in a naked singularity, was
presented in \cite{LehnerPretorius2010}. Since no fine-tuning of the initial data was
required, this result constituted a violation of the weak cosmic censorship, albeit
in spacetimes with compact extra dimensions. The
black rings of \cite{EmparanReall2002} also suffer from various types of instabilities
\cite{BRI}, including the Gregory--Laflamme  instability.

At very large angular momenta, black holes become highly
deformed and resemble black branes. The
rapidly
spinning Myers-Perry (higher-dimensional analogues of Kerr) black holes \cite{MP} in 
spacetime dimensions greater than 6 were shown to be unstable
under a (Gregory--Laflamme type of) 
``ultraspinning instability'' \cite{EmparanMyers}. 
In particular,
the end point of the axisymmetric ultraspinning instability of asymptotically flat Myers-Perry
6D black holes was studied in \cite{Figueras}, and this instability was 
found to give rise to a sequence
of concentric rings connected by segments of black membrane on the rotation plane which become 
thinner
over time in the non-linear regime, resulting in the formation of a naked singularity in finite time 
and 
consequently a violation of the
weak cosmic censorship conjecture.

\newpage

\subsection{Singularity resolution in GR by quantum effects} 

Einstein's theory suffers from the problem of
classical singularities, which are a generic feature of spacetimes in GR. The
existence of singularities indicate a breakdown of the classical theory at 
sufficiently large spacetime curvature, which is precisely when gravitational quantum effects are 
expected to be relevant. 
Consequently, QG is necessary for the clarification of whether the singularity theorems survive
when entering the quantum regime. 
The question
of whether a theory of QG can extend solutions of classical GR
beyond the singularities was first discussed in \cite{DeWitt}.

In any analysis of the singularity theorems in
the quantum realm, an important step
is the weakening of the energy conditions and finding an appropriate version of the curvature conditions. 
For example, averaged energy conditions to deal with the
quantum violations of the energy conditions have been considered.
It is also necessary to go beyond semiclassical theories and take into
account the quantum 
fluctuations of the spacetime itself, which leads to additional difficulties in seeking
quantum singularity theorems. In particular, in the classical theorems  the
pointwise focusing of geodesics is utilised, which cannot hold exactly (despite the smallness
of the  fluctuations) in a quantum regime. The notion of closed trapped surface can
also be generalized and adapted to quantum situations \cite{Wall2013}.

Let us discuss cosmological and black hole singularity resolution within loop
quantum gravity (LQG) and string theory.
LQG is a non-perturbative canonical
quantization of gravity based on Ashtekar variables \cite{Ashtekar2004},
in which classical differential
geometry of GR is replaced by a quantum geometry at the Planck scale,
and has been used to perform a rigorous quantization
for spacetimes
with symmetries.
Applying the techniques of LQG to cosmological
spacetimes is known as loop quantum cosmology (LQC),
in which the spatial homogeneity reduces the
infinite number of degrees of freedom to a finite number. 
Indeed, LQG reveals that singularities may be
generically resolved because of the quantum gravitational effects \cite{Singh14}. And due to the quantum
geometry, which replaces the classical differential geometry at the Planck scale, the
big bang is replaced by a big bounce, when energy density reaches a maximum value of about a half of the Planck density. 
The existence of a  viable non-singular bounce in the very early universe appears to be a generic result in all investigations of simple models of LQC, and
occurs without any violation of
the energy conditions or fine tuning \cite{AshtekarPawlowski}.

%%\newpage 

Often a singularity in GR,
such as the big bang and big crunch (to the future for contracting models)
as well as black holes, is characterised
by the divergence of a physical or geometrical
quantity (such as a curvature invariant)
and the
breakdown of the geodesic evolution. However, 
singularities can also arise due to pathologies of the tangent bundle, for
instance in conical singularities, or where there are directional singularities, defined as limit points towards
which the curvature tensor diverges along some (but not necessarily all) directions.
These complications led to an elaborate classification of possible singularities arising
from the curvature tensors \cite{EllisSchmidt}.
LQC attempts to resolve all singularities, including, for example,  the big rip,  and sudden
and big freeze singularities.

In contrast to the classical theory where singularities are a generic feature, there is 
growing evidence in LQC that singularities may be absent.
Recently various spatially homogeneous cosmological models have been studied within the
context of LQC \cite{AshtekarandSingh}. In particular, for the models that have been
exactly solved at the quantum level, the dynamics of sharply peaked, semiclassical states is
very well described by an effective theory that incorporates the main quantum corrections
to the dynamics \cite{Ashtekar,Diener2014}. For example, at the effective level
an infinite number of bounces and recollapses occur in the positive spatial
curvature FLRW model. 
Simple Bianchi type IX models have also been shown to be non-singular \cite{SinghWilson,CorichiMontoy}.
The original study \cite{Bojowald} was improved 
(to solve the problems with the infrared limit \cite{Ashtekar}), and the
Bianchi IX dynamics was constructed \cite{Ewing2010}. 
It is within the framework of the
improved dynamics that solutions to the effective equations for the Bianchi IX class of
spacetimes whose matter content is a massless scalar field was studied numerically  within LQC in \cite{CorichiMontoy}, and the big bang singularity was resolved and the
classical dynamics far from the bounce was reproduced.

%%\newpage

\vspace{3pt}

\noindent
{\bf{Problem P11:}} Can the singularity resolution results obtained in the spatially homogeneous spacetimes
in LQG be proven in a more general setting. 

\vspace{3pt}

Symmetry reduction within LQC entails a drastic simplification, and therefore 
important features of the theory might be lost by restricting the symmetry prior to quantization. 
However, it is believed that such studies do lead
to valuable hints on loop quantization  and inhomogeneous
spacetimes (and black holes) \cite{AshtekarandSingh}.
There is strong evidence from the numerical studies of
the BKL conjecture (see later), that near
the singularities the structure of the spacetime is not determined by the spatial derivatives,
so that it might be hoped that singularity resolution in spatially homogeneous models would
capture some aspects of the singularity resolution in more general  inhomogeneous spacetimes.
On the other hand, however, sometimes the limitations of LQC have been used to shed doubts on its
results. For example, in \cite{BojowaldPaily} it is claimed that a fully covariant
approach with validity beyond symmetry reduced scenarios produces physical results inequivalent
to those obtained from LQC (see also \cite{Singh14}). 
Recently LQG techniques have been used to study the effects
of QG in the simple Gowdy inhomogeneous  models 
with infinite degrees of freedom \cite{Brizuela}.
And the first steps in the study
of classical oscillatory singularities governed by the BKL dynamics
using LQG have been taken \cite{Garfinkle16}.

Loop quantization of black hole spacetimes uses similar techniques as in LQC, and leads
to similar results on singularity resolution 
\cite{AshtekarBojowald2006,GambiniPullin2008,Modesto,Corichi16}.
The resolution of gravitational black holes singularities  has been also studied
in string theory \cite{Risi,Cailleteau}. Indeed,
there has been significant progress on understanding
black holes in string theory recently  \cite{Lopes}, and some interesting examples have
been presented where gravitational singularities are resolved by higher derivative
corrections to the action \cite{Hub}. 
For example, the fundamental string in
five dimensions, which is singular in the standard supergravity description, is
regular after taking into account higher derivative corrections determined by
anomalies and supersymmetry \cite{Larsen}.
In particular, singularities were resolved in
string solutions of five dimensional supergravity corrected by the mixed gauge-gravitational
Chern-Simons term
with $AdS_3 \times S^2$ near string geometry 
(which can be interpreted in M-theory as M5-branes wrapped on four-cycles in a Calabi-Yau manifold)  \cite{Larsen}.
The techniques to resolve
singularities can applied in more general situations,
including black holes in five dimensions with different near horizon geometry,
rotating black holes, and generalizations
to other dimensions, including 10 and 11 dimensions, and theories with more general
matter content.

\newpage

\subsection{Black hole information paradox}

Hawking discovered that black holes are not completely black but emit a dim radiation due to quantum effects near the event horizon \cite{Hawking}. This result poses a 
fundamental theoretical problem because it appears to suggest that evaporating black holes destroy information, which is in conflict with a basic postulate of quantum mechanics that physical systems evolve in time according to the Schrodinger equation
(which is fully deterministic and unitary and thus no information can be truly lost or destroyed).
The apparent contradiction between Hawking's result and the unitarity postulate of quantum mechanics 
has become known as the {\em{black hole information paradox}} \cite{Hawking2005}
(see AL2.4 and AL2.5).

In more detail, although a black hole, formed by the gravitational collapse of a body in GR,
is classically stable, quantum particle creation processes will
result in the emission of Hawking radiation  \cite{Hawking} to infinity and corresponding mass loss of the black
hole, eventually resulting in the complete evaporation of the black hole. Semi--classical arguments,
from applying the local evolutionary laws of quantum field theory in a classical curved spacetime,
strongly suggest that in the process of black hole formation and evaporation, a pure quantum state
will evolve to a mixed state \cite{Marolf}. That is, 
if the black
hole itself has completely disappeared then only the thermal radiation is left,
and this final state would be largely independent of the initial state and
would thus not suffice to deduce the initial state and information would have been lost.
But this behavior is quite different from that of familiar quantum systems
under unitary time evolution. 
There are, however, a number of natural ways to attempt to restore unitarity, including
QG corrections and additional degrees of freedom, in addition to a modification of quantum mechanics itself.

\newpage

The black hole information paradox is really a combination of two problems: the causality paradox and the entanglement problem.
For illustration, consider a spherical shell of mass collapsing to form
a black hole. In the semiclassical approximation the shell passes through its horizon,
and ends at a singularity. 
(A) After the shell passes
through its horizon, light cones in the region between the
shell and the horizon ‘point inwards’. If we assume that ‘faster
than light’ propagation is not possible, then the information
in the shell is causally trapped inside the horizon.
Thus this information cannot escape to infinity as the
hole evaporates away. 
(B) The process of Hawking
radiation creates entangled pairs at the horizon.
But the large entanglement between
the radiation near infinity and the remaining hole
near the endpoint of evaporation
may be beyond the semi--classical approximation.

The various proposals to resolve the information paradox include the idea that
quantum fields near black hole horizons are not in
fact well-described by vacuum but are instead highly excited due to new unknown
physics. The set of excitations is called a black hole ``firewall," and might even be sufficiently strong that
spacetime fails to exist in any recognizable sense in the interior of such black holes.
This has also been discussed within the context of gauge/gravity duality.
In particular, the AdS/CFT correspondence \cite{Klebanov} partially resolves the black hole information paradox 
since it describes how a black hole 
(i.e., particles on the boundary of AdS spacetime)
can evolve in a unitary fashion in a manner consistent with quantum mechanics,
leading to information conservation  in this context (see also \cite{Strominger}).

Other alternatives to information loss include:
(i) A black hole never actually
forms in the collapse, but rather some
other structure without an event horizon, such as a ``fuzzball'', is formed.
(ii) A black hole forms in the expected manner but there are major
departures from semi--classical theory and there is greatly diminished entanglement during the evaporation process.
(iii) The evaporation process shuts off by the time the black hole has evaporated down to the
Planck scale when QG effects become dominant, such that the resulting “remnant”
contains all of the “information” that went into the black hole.
(iv) The evaporation
process proceeds as in the semiclassical analysis until the black hole reaches the Planck scale, whence
all of the information that had been stored within the black hole then
emerges in a final burst.

In another older but still plausible suggestion  \cite{Susskind93}
a phenomenological description of
black holes and their event horizons  was introduced
based upon three postulates which, when  implemented
in a ``stretched horizon''  description of a black hole, preserve
free infall ``through the
horizon" within the full quantum
theory. 
The notion of black hole ``complementarity'' is then realised,
whereby an observer outside of the black hole receives the
information returned from the horizon to infinity (in the form of radiation emanating from the apparent horizon,
which is presumed to be outside the event horizon for a dynamical shrinking black hole),
but observers inside the black hole cannot communicate
with the outside. Therefore, any possible contradictions 
might be acceptable since they are not visible to any single distant observer and consequently
there would be
no  resulting tension with any known experiments.

\newpage

\section{Problems in Cosmology (Ph5)} 

Cosmology is the study of the large scale behaviour of the Universe within a theory of gravity, 
which is usually taken to be GR.
There are many open problems in theoretical cosmology.
For example, what precisely is
the hypothetical inflaton field and
what are the details of 
{\em{cosmic inflation}}?  Is inflation self-sustaining
through the amplification of quantum-mechanical fluctuations and thus still occurring in some 
(distant) places in the Universe?
Does it give rise to countless ``bubble universes'' and, if so,
under what initial conditions, and   does a {\em{multiverse}} exist?
Cosmological inflation is generally accepted as a solution to the {\em{horizon problem}},
that the universe appears more uniform on larger scales than expected,
but are other  explanations  possible?
What is the origin and future of the
universe and, in particular,  is the universe heading towards some sort of final singularity?
Or is it evolving towards a big bounce or is it even part of an
infinitely recurring cyclic model?

Since cosmology concerns the behaviour of the Universe when the small-scale structures such as stars and galaxies
can be neglected, the ``Cosmological Principle'' (a generalization of the 
Copernican principle)  is often assumed to hold, which asserts that: 
{\em{On large scales the Universe can be well--modeled by a solution to Einstein's equations which is spatially homogeneous and isotropic.}} That is, 
a (possibly preferred) notion of cosmological time can be picked such that at every instant  on large scales
space looks identical in all directions (isotropy), and (spatial homogeneity) it is not possible to distinguish between any two points (which is clearly not true on the astrophysical scales of galaxies).
However, it would be more satisfactory if the cosmological principle could be derived as a consequence of GR 
(under suitable assumptions), rather than something assumed a priori. 
That is, could spatial homogenization and isotropization at late times  be derived as a mathematical consequence of Einstein's equations under appropriate physical conditions and for suitable initial  data.
This question is partially addressed within the inflationary paradigm.

\paragraph{Dark matter and dark energy:}

Perhaps the most important questions in cosmology are those concerning dark matter and dark energy.
These types of matter and energy are detected by their gravitational influences, but can't be observed directly.  
The estimated distribution of dark matter in the Universe is based on observed galaxy rotation curves,
nucleosynthesis predictions and structure formation computations \cite{Freese}. Although the
identity of the missing  dark matter
is not yet known (e.g., whether it is a particle, perhaps the lightest superpartner, or
whether the phenomena attributed to dark
matter is not  described by some form of matter but rather by an extension of GR), it is generally believed that this problem will be solved by conventional physics. The dark energy problem is much more serious.
Indeed, this problem is widely regarded as one of the major
obstacles to further progress in fundamental physics
\cite{Witten2001,Steinhardt}.

The {\em{cosmological constant  problem}} was  discussed comprehensively by
Steven Weinberg \cite{Weinberg1989}. Standard quantum field theory
(QFT) predicts a huge vacuum energy density from
various sources. But the equivalence principle of
GR requires that every form of energy gravitates in the
same way, so that it
is widely believed that the vacuum energy gravitates
as a cosmological constant
which would then have an
enormous effect on the curvature of spacetime.
However, the observed effective
cosmological constant is so small compared with
the  predictions of QFT that an unknown bare cosmological
constant has to cancel this huge contribution from the
vacuum to better than up to at least 120 decimal places (AL2.1).
It is an extremely difficult fine-tuning problem that gets
even worse when the higher loop corrections are included
\cite{Padilla}.
More recently Weinberg and others have adopted the view
that, of all of the proposed solutions to this problem, the only acceptable one is the controversial
anthropic bound  \cite{Weinberg1987}.

In addition, the  Universe has been accelerating
in its expansion for the last few billion years \cite{Riess,Perlmutter}. 
Within standard cosmology  the cause of this acceleration is commonly called {\em{dark energy}}, which appears to
have the same properties as a relatively tiny cosmological constant, an
effectively repulsive gravitational force (or {\em{levitational force}}) in GR. 
The additional
cosmological constant (coincidence) problem of explaining why it has
such a specific small observed value, which is
the same order of magnitude as the present mass density
of the matter in the Universe, must also be faced (AL2.2).
It is often speculated as to whether
dark energy is a pure cosmological constant or whether dynamical models such as,
for example, quintessence and phantom
energy models are more appropriate. 
Some physicists have also proposed alternative
explanations for these gravitational influences, which do not require new forms of matter and energy, but these
alternatives are not popular and lead to modified gravity on large scales.
The possible cause of the observed acceleration of the Universe has also been discussed within the context of backreaction and inhomogeneities (see later). 

Finally, it has also been proposed that an
observed ``dark flow'',  a non-spherically symmetric gravitational pull from outside the observable Universe, is responsible
for some of the observed motion of large objects such as galactic clusters in the Universe.
Analyses of the local bulk flow of
galaxies (as measured in the frame of the cosmic microwave background) indicate a lack of 
convergence to the cosmic background frame
even beyond 100 Mpc \cite{Kashlinsky}, in contrast to standard
expectations if the Universe is in fact spatially
homogeneous on larger scales.
Indeed, low redshift supernova data have shown that
there is an anomalously high and apparently constant bulk flow of approximately 250 km/s extending
all the way out to the Shapely supercluster at approximately 260 Mpc and further, a 
discrepancy which has been confirmed by
analysis of the 6dF galaxy redshift survey 
\cite{Sarkar}.

%%\newpage

However, although mathematics is very important in many of these problems, they are not problems in mathematical physics per se.
In addition, numerical computations have always played an important role in physical cosmology \cite{Computational},
but it is not clear that such calculations are within the remit of mathematical physics. Indeed,
computational cosmology
within full GR is now beginning to
address fundamental issues \cite{Bentivegna,Giblin,Adamek}. Let us briefly discuss some
topical examples.

Studies of ``bubble universes" in
which our Universe is one of many, nucleating and growing inside an
ever-expanding false vacuum, have been undertaken with computational cosmological tools.  In
particular, \cite{bubbles} investigated the collisions between bubbles.
It is expected that initial conditions will contain some measure of inhomogeneities,
and random initial conditions will not necessarily give rise to
an inflationary spacetime. 
It has been shown 
that large field inflation is robust to simple inhomogeneous (and
anisotropic) initial conditions with large initial gradient energies in
situations in which the field is initially confined to the part of the
potential that supports inflation, while 
it is also known that small
field inflation is much less robust to inhomogeneities than its large field counterpart \cite{infl}.

Using exotic matter, or alternative modified theories of gravity,
can classically lead to the initial singularity being replaced
by a {\em{bounce}} to an expanding universe \cite{Brandenberger}.
For example, computational cosmology
methods have been applied to the study of bouncing cosmologies
in the ekpyrotic cosmological
scenario; by studying
the evolution of adiabatic perturbations in a nonsingular bounce \cite{EKp}, 
it was  shown that the bounce is disrupted in regions 
with significant spatial inhomogeneity and anisotropy compared with the background
energy density, but is achieved in regions that are relatively spatially homogeneous and
isotropic.  The precise properties of a cosmic bounce depend upon the way in which it is
generated, and many mechanisms have been proposed for this both classically and non-classically. 
There are
quantum gravitational effects associated with string theory \cite{Turok} and LQG \cite{Bojowald,Ashtekar}.
In particular, in LQC there is
such a bounce when the energy density reaches a maximum value of
approximately one half of the Planck density.

%%\newpage

Nevertheless, some precise mathematical questions  in cosmology can be formulated.
For example, there are questions about the generality of inflation for generic initial data
(although precise statements are difficult because there are many theories of inflation and
there are no natural initial conditions).
But mathematical  theorems are  possible in the study of the stability of de Sitter spacetime.
This is part of the more general question of the {\em{stability of cosmological solutions}}
(namely,  if 
a cosmological solution is perturbed a little bit by, for example, factoring in the small-scale structure, 
is the evolution as governed by Einstein's 
equations qualitatively the same in the large as the evolution of the underlying cosmological solution). 
This requires the study of the (late time)  behaviour of a complicated set of partial differential equations around a special solution (and there are several cosmological models that are of particular interest, including the very
simple Milne model \cite{AnderssonMoncrief,LAndersson14}). 
These are genuinely problems in mathematical physics.

We first recall  that when
the cosmological constant vanishes and the matter satisfies the usual energy
conditions, spacetimes of Bianchi type IX recollapse and so are never
indefinitely expanding.
This is formalized in the
“closed universe recollapse conjecture” \cite{BarrowTipler}, which
was proven by Lin and Wald  \cite{LinWald}.
However,  spacetimes of Bianchi type IX need not always recollapse when there is a non-zero positive 
cosmological constant.

\newpage
\subsection{Stability of de Sitter spacetime}

In \cite{Friedrich1986} Friedrich proved (using regular conformal field equations) a result on the stability 
of de Sitter spacetime: all initial data (vacuum with positive cosmological constant) near
enough (in a suitable Sobolev topology) to initial data
induced by de Sitter spacetime on a regular Cauchy hypersurface
have maximal Cauchy developments
which are geodesically complete.
de Sitter spacetime is thus an attractor for expanding cosmological
models with a positive cosmological constant. The result also gives additional details
on the asymptotic behaviour and may be thought of as proving
a form of the so-called `` cosmic no hair'' conjecture in the vacuum case.
For more recent work see \cite{Ringstrom2015} and references within.

A general theorem of Wald \cite{Wald83} states that any spatially homogeneous model
whose matter content satisfies the strong and dominant energy conditions and
which expands for an infinite proper time (i.e., does not  recollapse) is asymptotic to an isotropic de Sitter spacetime. This 
cosmic no hair theorem
does not depend on the details of
the matter fields, and therefore the question remains as to whether solutions corresponding
to initial data for the Einstein equations with a positive cosmological constant
coupled to reasonable matter exist globally in time only under the condition that
the model is originally expanding. It can be shown that this is true for various
matter models using the techniques of \cite{Rendall95}.

Models with a scalar
field with an exponential potential are also inflationary because the rate of (volume)
expansion is increasing with time, and global results are possible \cite{Coleybook,exppot}.
Inflationary
behaviour also arises in the presence of a scalar field with a power law
potential, but occurs at intermediate times rather than at late times.
Local results are then possible but are difficult; primarily this question is studied numerically.
 
It is of interest to know what happens to the cosmic no--hair theorem 
in inhomogeneous spacetimes. 
Some partial results are possible for a positive cosmological constant in the inhomogeneous case \cite
{Jensen}.
But even less is known for scalar field models with an exponential potential  \cite{Coleybook}.
\vspace{3pt}

\noindent

{\bf{Problem P12:}} Prove a cosmic no--hair theorem 
in generic inhomogeneous spacetimes.
 
\vspace{3pt}
 
%%\newpage 

The potential instability of de Sitter spacetime in quantized theories has been investigated. 
In a semi-classical analysis of backreaction in an expanding Universe with a conformally coupled scalar 
field and vacuum energy it was shown that a local observer perceives de Sitter spacetime to 
contain a constant thermal energy density despite the dilution from expansion due to a continuous 
flux of energy radiated from the horizon, leading to  the  evolution of the Hubble rate at late times
which deviates significantly from that in de Sitter spacetime, which is thus
unstable \cite{Markkanen}. 
This result is in apparent disagreement with the
thermodynamic arguments in \cite{GibbonsHawking77} in which it was concluded that unlike black holes de
Sitter spacetime is stable. However, 
if de Sitter spacetime were unstable
to quantum corrections and could indeed decay, it could provide an important mechanism
for alleviating the cosmological constant problem and perhaps also the fine-tuning issues
encountered in the extremely 
flat inflationary potentials that are required by observations.
A de Sitter instability would certainly have a profound impact on the fate of the
Universe since it rules out the possibility of an eternally exponentially expanding de Sitter
spacetime as classically implied by the standard concordance model.
This issue is currently unresolved.

\newpage

\subsection{Cosmological singularities and spikes}

The singularity theorems tell us that singularities occur under very general circumstances in 
GR, but they say little about their nature  \cite{Senovilla2012}.
Belinskii, Khalatnikov and Lifshitz (BKL) \cite{art:LK63} 
have conjectured that within GR, for a generic inhomogeneous cosmology, the approach to the
(past) spacelike singularity is vacuum dominated, local, and oscillatory (mixmaster).
The associated dynamics 
is referred to as the BKL  dynamics. In particular, 
due to the nonlinearity of the
Einstein field equations, if the matter is not a  (massless) scalar field,
then sufficiently close to the singularity one can {\em{neglect all matter terms}} in
the field equations relative to the dynamical anisotropy.
BKL checked
that their assumptions were consistent with the Einstein field equations; but
that doesn't necessarily mean that those assumptions hold in general physical situations.
Recent 
numerical simulations have verified that the BKL conjecture is satisfied for 
special classes of spacetimes \cite{Berger,DavidG}.

To date there have been three main approaches to investigate the structure of generic singularities,
including the heuristic BKL metric approach and the Hamiltonian approach.
The dynamical
systems approach \cite{WE}, in which 
Einstein's field equations are recast into scale invariant asymptotically regularized dynamical systems (first order systems of autonomous ordinary differential
equations and partial differential equations) in the approach towards a generic spacelike singularity, offers a
more mathematically
rigorous  approach to cosmological singularities.
In particular, a dynamical systems formulation for the Einstein field equations without
any symmetries was introduced in \cite{Uggla03}, resulting in
a detailed description of the generic attractor, concisely formulated conjectures
about the asymptotic dynamic behavior toward a generic spacelike singularity, and a basis
for a numerical investigation of generic singularities \cite{Andersson}.

In more detail, in order to construct the solution in a sufficiently small spacetime neighborhood
of a generic spacelike singularity \cite{Heinzle,Uggla03} Einstein's field equations
are reformulated by
assuming that a small neighborhood
near the singularity can be foliated with a family of spacelike surfaces such that
the singularity ``occurs'' simultaneously, and
the expansion of the
normal congruence to the assumed foliation  
are factored out by utilizing a conformal transformation
(whence the Einstein field equations split
into decoupled equations for the conformal factor and a coupled system of dimensionless
equations for quantities associated with the dimensionless conformal metric).

\newpage

Unfortunately, until recently very few rigorous mathematical statements had been made. 
Based on the work of Rendall \cite{Rendall},  Ringstrom produced the first major proofs about asymptotic  spatially homogeneous Bianchi type IX
cosmological dynamics \cite{Ringstrom}. 
Notably, Ringstrom obtained the first theorems about
oscillatory behavior of generic singularities for Bianchi type VIII and, more substantially, type IX models in GR. In particular, Ringstrom managed to prove that the past attractor
in Bianchi type IX resides on a subset that consists of the union of the Bianchi
type I and II vacuum subsets.
But this theorem does not identify
the attractor completely, nor determine if the Kasner map is relevant for dynamics asymptotic
to the initial singularity in Bianchi type IX, and the theorem says very little
about Bianchi type VIII \cite{Heinzle}
(however, see the recent work of Brehm \cite{Brehm}).

In the spatially homogeneous case the focus is on mathematically rigorous results. For example, it has been argued that the idea that Bianchi type IX models are essentially understood is a misconception,
and what has actually been proven about type IX asymptotic dynamics was addressed in  \cite{Heinzle}
(however, see \cite{Brehm}).
In particular, all claims about chaos in Einstein's equations (especially at a generic spacelike
singularity) rest
on the (plausible) belief that the Kasner map 
(which is associated with chaotic properties)
actually describes the asymptotic dynamics of
Einstein's equations; however,
this has not been proved rigorously to date \cite{Heinzle}. Most importantly,  the
role of type IX models in the context of generic singularities not yet been rigorously established                          \cite{Heinzle09,Uggla03,Heinzle}.
There remain several important open problems, including:

\vspace{3pt}

\noindent
{\bf{Problem P13:}}  Prove that the past attractor of the Bianchi type IX dynamical
system coincides with the Mixmaster attractor (as defined in \cite{Heinzle} -- the Bianchi II variety of \cite{Ringstrom}) rather than being a subset thereof. 

\vspace{3pt}

The BKL oscillatory dynamics have been studied in simple perfect fluid models with a linear equation of state.
Some matter fields can have an important effect on the dynamics near the
singularity. A scalar field or stiff fluid leads to the oscillatory behaviour being
replaced by monotonic behaviour and consequently
to a significant reduction in the complexity of the dynamics  close to the singularity \cite{BGZK}. Based on
numerical work and the qualitative analysis of \cite{Hewitt2003},
the so–-called exceptional Bianchi type VI$_{-\frac{1}{9}}$ class B model
(which  has the same number
of degrees of freedom as the most general  Bianchi type VIII and
IX class A models) has an oscillatory
singularity. An electromagnetic field can
lead to oscillatory behaviour which is not present in vacuum or
perfect fluid  models of the same symmetry type. For example, models of Bianchi
types I and VI$_0$ with an electromagnetic field have  oscillatory behaviour
\cite{LeBlanc}. Oscillations can also occur in all Bianchi models in the presence of a tilting fluid 
\cite{Hewitt2001,Hervik07}.

\newpage

It is imperative to
discuss
generic oscillatory singularities in inhomogeneous cosmologies.
In \cite{Uggla13}
qualitative and numerical support was presented for the BKL scenario in the
Hubble-normalized state space context for an open set of time lines. 
In more generality, a  heuristic physical justification of asymptotic locality may
be that ultra strong gravity causes
particle horizons shrink to zero size toward the singularity along each timeline, which
prohibits communication between different time lines in the asymptotic limit (and may
hence be referred to as asymptotic silence).

To gain further insights about general oscillatory singularities in inhomogeneous
spacetimes, models with two
commuting spacelike Killing vectors (so-called $G_2$ models) have been investigated. 
The BKL dynamics has been discussed in
generic vacuum, spatially compact $U(1) \times U(1)$–-symmetric
spacetimes with vanishing twist  and in generic polarized $U(1)$ spacetimes 
\cite{LARS99}, and in twisting $U(1) \times U(1)$--symmetric vacuum models on $T^3$
Gowdy models, and on $S^2 \times S^1$, $S^3$
and lens spaces L(p, q) \cite{Maier}.

The description of generic asymptotic dynamics towards a
generic spacelike singularity in terms of an attractor, has resulted in mathematically precise
conjectures \cite{Heinzle09,Uggla03}, and involves the possible existence of
finite dimensional attractors in infinite dimensional systems \cite{Temam}.

\vspace{3pt}

\noindent
{\bf{Problem P14:}}
Prove the BKL  locality
conjecture in the general inhomogeneous context.

\vspace{3pt}

\newpage

\paragraph{Spikes:}

Recently, a new spike phenomenon that had not been anticipated by BKL
was found in numerical simulations  \cite{Berger}.
Since it is a general feature of solutions of partial differential equations that spikes
occur it is, of course, expected that they  occur in solutions of Einstein's field
equations in GR. 
In the case of spikes, the spatial derivative terms do have a significant
effect at special points. 
In particular, in the approach to the singularity in the mixmaster regime,
a spike occurs when a particular
spatial point is stuck in an old Kasner epoch while its neighbours eventually bounce
to the new one.  
Because spikes become arbitrarily
narrow as the singularity is approached, they are a challenge to the numerical
simulations. Spikes are also a challenge to the mathematical treatment of 
spacetimes. Mathematical justification has been presented in  \cite{Ringstrom2004}.  
More success has been obtained in finding 
exact spike solutions ~\cite{art:Lim2015}.

Spikes were originally found numerically in the context of vacuum orthogonally transitive,
spatially inhomogeneous $G_2$ models
~\cite{Berger,art:Limetal2009}. 
Therefore,  numerical studies of $G_2$ and more general cosmological models have produced evidence that the BKL conjecture generally holds except possibly at isolated points 
(surfaces in the three-dimensional space) where spiky structures (``spikes'') form 
\cite{art:Bergeretal2001}, and the asymptotic locality part of the BKL conjecture is violated.
Spikes naturally occur in a class of non-vacuum $G_2$ models and, due to gravitational instability, leave
small residual imprints on matter in the form of matter perturbations.  Particular interest has been paid to
spikes formed in the initial oscillatory regime, and their imprint on matter and
structure formation has been studied  numerically \cite{art:ColeyLim2012}.  

Therefore, 
generic singularities are not only associated with asymptotic locality but also with non-local recurring spikes, although it is believed that a set of measure zero of
timelines exhibit spike formation \cite{Andersson}.

\vspace{3pt}

\noindent
{\bf{Problem P15:}} Prove the existence of spikes and determine their effect on any eventual
generic singularity proofs.
\vspace{3pt}

\newpage

There are other unresolved questions pertaining to recurring spike behavior and generic spacelike singularities. For example, 
are there spikes that undergo infinitely many recurring spike transitions? How, where and how often
do spikes form? How does spike interference work and
can spikes annihilate? Are there generic singularities without recurring spikes and are
there generic singularities with a dense set of recurring spikes? 
Some of these issues have been discussed recently in \cite{HeinzleUggla2012}.

The asymptotic dynamics of general
solutions of the Einstein vacuum equations toward a generic spacelike
singularity have been studied. Matter sources such as spatially homogeneous  perfect fluids and simple massless scalar fields
\cite{art:LK63,WE} have been considered. 
It is of
particular interest to determine the structural stability of generic inhomogeneous spacelike singularities, 
especially by including matter such as
massless scalar fields
and electromagnetic fields (which influence the generic spacelike singularity in different ways),
and to also go beyond
GR and include form fields.

\vspace{3pt}

\noindent
{\bf{Problem P16:}} Determine the structural stability of generic inhomogeneous spacelike singularities 
for general matter fields present in the early universe.
\vspace{3pt}

In  \cite{Damour}, a heuristic Hamiltonian approach 
(closely
connected to that of BKL) was used to study the dynamics
of the Einstein-dilaton-p-form system in the neighborhood of a generic spacelike
singularity. The asymptotic behavior of the fields was described  by a ‘billiard’
motion in a region of hyperbolic space bounded by straight ‘walls’ (dubbed ``cosmological billiards''), 
and a remarkable connection
between the asymptotic dynamics of generic spacelike singularities and Kac-Moody algebras was revealed \cite{Damour}.
A link between the Hamiltonian
and the dynamical systems approach
to inhomogeneous cosmologies was established in  \cite{Heinzle09}.
More recently the fermionic sector of supergravity theories, in which the  gravitino is treated
classically, was studied \cite{Kleinschmidt}. The  quantum
generalization of the resulting  fermionic cosmological billiards,
defined by the dynamics of a quantized supersymmetric
Bianchi type IX cosmological model (within simple 4D supergravity)
\cite{DamourSpindel}, was also investigated. The
hidden Kac-Moody structures were again displayed. 

\newpage

\paragraph{Isotropic singularity:} Based on entropy considerations, 
Penrose  \cite{Penrosebooks} proposed the ``Weyl curvature hypothesis'' that asserts that the initial singularity
in a cosmological model should be such that the Weyl curvature tensor tends to
zero or at least remains bounded. There is some difficulty in representing this condition 
mathematically and it was proposed in \cite{Goode} that a clearly formulated
geometric condition, which on an intuitive level is closely related to the original
condition, is that the conformal structure should remain regular at the singularity.
Singularities of this type are known as conformal or isotropic singularities.
It has been shown \cite{Claudel,Newman} that
solutions of the Einstein equations coupled to a perfect fluid satisfying the
radiation equation of state with an isotropic singularity are determined uniquely by certain free
data given at the singularity. The data which can be given is, roughly speaking,
half as much as in the case of a regular Cauchy hypersurface.
In \cite{Anguige} this was extended to the linear equation of state
case, and can be extended to  more general matter (e.g., general fluids and a collisionless gas of 
massless particles) \cite{Rendall2002}.

Many additional questions can be asked in the context of alternative, modified theories of gravity such 
as, for example, the general applicability of the BKL behaviour close to the cosmological 
singularity.
Such questions will not be included here. However, the following mathematical physics question
on isotropization is relevant.

\vspace{3pt}

\noindent
{\bf{Problem P17:}} Are  isotropic singularities typical in modified theories of gravity.

\vspace{3pt}

The stability of the isotropic vacuum Friedmann universe  on 
approach to an initial cosmological singularity in gravity theories with higher--order 
curvature terms added to the Einstein-Hilbert Lagrangian of GR have been studied \cite{Middleton}. A special
isotropic vacuum solution exists which behaves like the radiation-dominated Friedmann universe and is
stable to anisotropic and small inhomogeneous perturbations in the past, unlike the situation in GR.
An analytical solution valid for particular values of the equation of state parameter was also found 
such that the singularity is isotropic in a higher dimensional flat anisotropic Universe 
filled by a perfect fluid in Gauss-Bonnet gravity 
\cite{Kirnos}. Some simple cosmological solutions of gravity theories 
with quadratic Ricci curvature terms added to the Einstein-Hilbert Lagrangian  
have also be studied
\cite{BarrowHervik}.

\newpage

\subsection{Averaging Einstein's field equations}

The averaging problem in GR is of fundamental importance.
The gravitational field equations on large scales are obtained by averaging
or coarse graining the Einstein field equations of GR. 
The averaging problem in cosmology is crucial for the correct
interpretation of cosmological data. 
The so--called fitting problem is perhaps the most important unsolved
problem
in mathematical cosmology \cite{fit}.

The spacetime or space volume averaging approach  must be well posed and generally covariant.
This raises important new questions in differential geometry. 
The formal mathematical issues of averaging tensors on a differential manifold have recently
been revisited \cite{Av,Coley10,Averaging}. 
The coarse grained or averaged field equations need not take on the same mathematical
form as the original field equations. Indeed, in the
case of the macroscopic gravity approach \cite{Averaging,CPZ} the
averaged spacetime is not necessarily even Riemannian.
Scalar quantities can be averaged in a straightforward manner.
In general, a spacetime is completely characterized
by its scalar curvature invariants, and this suggests a particular spacetime averaging
scheme based entirely on scalars \cite{Coley10}.
In the approach of
Buchert \cite{bu00} a 1+3 cosmological spacetime splitting is employed and only scalar quantities
are averaged.

The spacetime averaging
procedure adopted in macroscopic gravity, 
which is fully covariant and gauge independent,
is based on the concept
of Lie-dragging of averaging regions, and it has been shown to exist on an arbitrary Riemannian spacetime with well-defined local averaged properties
(however, see \cite{Av}). The averaging of the structure
equations for the geometry of GR then produces the structure equations for the averaged
(macroscopic) geometry and 
gives a prescription for the (additional) correlation functions (in the macroscopic field equations) which
emerge in an averaging of the non-linear field equations \cite{Averaging}.

\vspace{3pt}

\noindent
{\bf{Problem PF2}} Provide a rigorous mathematical definition for averaging in GR.
\vspace{3pt}

Although the standard spatially homogeneous and isotropic Friedmann--Lemaitre--Robertson--Walker (FLRW) model
(or so-called $\Lambda$CDM cosmology) has been extremely successful in describing current observations
(up to various possible
anomalies and tensions \cite{tension}, and particularly the existence of structures on gigaparsec
scales such as the cold spot and some super-voids \cite{Finelli}), it requires sources of dark energy density that dominate the present Universe that have never 
been directly detected.
In addition, the Universe is not isotropic
or spatially homogeneous on local scales. An averaging of inhomogeneous spacetimes on
large scales can lead to important effects. For example, on cosmological scales the dynamical
behavior can differ from that in the standard FLRW model and, in
particular, the expansion rate may be significantly affected \cite{bu00}.

%%%%%%%%%%

\newpage

Indeed, current observations of the structure of the late epoch Universe reveal a
significantly complex picture in which groups and clusters of galaxies of various sizes
form the largest gravitationally bound structures, which themselves
form knots, filaments and sheets that thread and surround very underdense
voids, creating a vast cosmic web \cite{web}. A significant fraction of the volume of the
present Universe is in voids of a single characteristic diameter of approximately $30$ megaparsecs
\cite{HV1} and a density contrast which is close to being empty, so that by volume the present universe is
void--dominated \cite{Pan11}.

A hierarchy of steps in coarse graining
is necessary
to model the observed  complex gravitationally bound  structures on large scales \cite{dust}.
In the standard FLRW cosmology it is implicitly assumed that regardless
of the gravitational physics in the coarse graining hierarchy,
at the final step  the matter distribution
can be approximated by an ``effective averaged out'' stress-energy tensor.
However, the smallest scale on which a notion of statistical
homogeneity arises is $70$--$120$ megaparsecs \cite{sdb12}, based on the two-point
galaxy correlation function, and  variations of the number density of
galaxies of order 7--8\% are still seen when sampling on the largest possible
survey volumes \cite{h05,sl09}.

\vspace{3pt}

\noindent
{\bf{Problem PF3}} Can averaging play an important role in cosmology. In particular,
what is the
largest scale that we can coarse-grain matter and geometry that obeys
Einstein's equations on smaller scales such that the average evolution
is still an exact solution of Einstein's equations.

\vspace{3pt}

After coarse graining we obtain a smoothed out macroscopic geometry (with
macroscopic metric) and macroscopic matter fields, valid on larger scales.
Indeed, the averaging of the Einstein field equations for local inhomogeneities on
small scales can in general lead to very significant dynamical {\em backreaction} effects
\cite{NotGW} on the average evolution of the
Universe \cite{bu00}.
In addition, all deductions about cosmology are based on light paths (null geodesics)
that traverse great distances (which  preferentially travel  through underdense
regions -- the voids in the real Universe). However,
inhomogeneities affect curved null geodesics and can significantly alter
observed luminosity distances. This leads to the following fundamental problem: although photons follow null geodesics in the local geometry, what trajectories
do photons follow in the averaged macro-geometry \cite{AAC2}. More importantly, however, is the fact that
averaging (and inhomogeneities in general)
can affect the interpretation of cosmological data \cite{BC,clocks,AAC2}.

A topical but theoretically conservative approach is to 
treat GR as a {\em mesoscopic} (classical) theory
applicable on those small scales on which it has actually
been tested, with a local metric field (the geometry) and matter fields, whence 
the effective dynamical equations on cosmological scales are then obtained by averaging.
In this approach, backreaction effects might
offer a resolution to problems related to dark energy and dark matter.

\newpage

%%%%%%%%%%%%%%%%%%%%%%%%%%%%%%%%%%

\section{Summary of problems in mathematical physics}

The hardest part, perhaps, is making the number of open questions in mathematical physics add up to 42. 
First, we must decide whether we mean 
42 in the mathematical sense (i.e., the exact number 42), or in the physics sense (i.e.,
a number between 40 and 44).

There are the 5 classical problems, 
H6, S3, S8, S15 and M2, and the related problems BS1 and BS3 (and the problem of turbulence)
and the more specific problems in YM theory and their generalizations to EYM theory (e.g., Y117 and Y118). 
In addition, there are the problems BS2, BS8 and BS14. Most of the problems proposed by Bartnik and Penrose
have been subsumed in the open problems P1 -- P17. However, problems RB20, RB21, RB32 (Y115) and RB43 remain.
There are 7, 4 and 6 open problems in each of GR, the quantum realm and cosmology, respectively, in the list P1 -- P17. I have also listed 3 personal favorites (PF1 -- PF3).

As mentioned earlier, 
many of the most important problems in theoretical physics are generally  not problems in mathematical
physics, despite the fact that a lot of mathematics is often utilized
(as discussed earlier, for example, in string theory).
Some of the problems which are absolutely fundamental for 
theoretical physics, and which almost by definition
are vague and not yet well formulated, have been briefly discussed in the text. 
But they
may (or may not)  turn into bone fide problems in mathematics or mathematical physics in the future ({\bf{MPF}}).
For example, although important, the question to explain the
anthropic reasons for the fine tuning of our Universe is not likely to to lead to an explicit problem in
mathematical physics. The following important problem may well lead to a
problem in mathematical physics.

\vspace{3pt}

\noindent
{\bf{Problem MPF1:}} Resolve the black hole information paradox. 

\vspace{3pt}

The two most fundamental problems in
theoretical physics will likely lead to
problems in mathematical physics in the future (see AL1/AL2 and AL3).

\vspace{3pt}

\noindent
{\bf{Problem MPF2:}} The cosmological constant problem and dark energy.
\vspace{3pt}

\vspace{3pt}

\noindent
{\bf{Problem MPF3:}} Formulate a fully consistent theory of QG.
\vspace{3pt}

\newpage

Numerical computations have always played an important role 
in any mature area of theoretical physics (such as GR and more recently 
in computational cosmology \cite{Computational}).
But it is not clear that, in general, such numerical problems are really problems in
mathematical physics. In addition,
numerical problems typically also require the
``complete control'' on the behaviour of gravitation in the very non-linear regime.
This always concerns the technical and practical question of the  ``cost'' due to the length of the required
computation and the small numerical error necessary to
ensure the solution can be trusted,  
which is not really a problem in  mathematical
physics.

There are perhaps no
important open questions in numerical relativity per se. 
On the other hand,
any important problem within GR that involves nonlinear phenomena 
would be an important problem for numerical
relativity.
Numerical work supports many of the conjectures discussed in this paper and has led to many important  theoretical
advances. For example, the mathematical stability of AdS spacetime has an important 
numerical component and
cosmic censorship is supported by numerical experimentation.
In addition, we have discussed the role of numerics in 
the understanding of spikes and the BKL dynamics, 
and in problems in cosmology and higher dimensional gravity.
In particular, we discussed the problems of the
generality of
bouncing models (at a cosmological singularity)
and of inflation for generic initial data. The latter problem may well lead to a genuine 
computationally motivated problem in mathematical
physics ({\bf{CMP}}), at least within a specific physical realization of inflation.

\vspace{3pt}

\noindent
{\bf{Problem CMP1:}} What is the generality of inflation for generic initial data. 
\vspace{3pt}

There are also important numerical problems in relativistic astrophysics,
such as in  the ultra-relativistic regime of interactions,
infinite-boost black hole collisions and colliding plane-fronted waves,
and most importantly in black hole mergers in general.
As noted earlier, the two-body problem has played and continues to play 
a pivotal role 
in gravitational physics \cite{Focusissue}. 
Recent advances in numerical computations have enabled the study of the 
violent
inspiral and merger of two compact objects (such as, for example, black 
holes and 
neutron stars), in which 
an enormous amount of gravitational radiation is produced.

In particular, the 
detection and analysis of the gravitational-wave signals generated 
by black hole collisions necessitate very precise
theoretical predictions for use as 
template waveforms to be cross-correlated against the output of 
the gravitational-wave detectors, which is of great importance in light 
of recent
LIGO observations \cite{LIGO}.
The orbital dynamics and gravitational-wave emission of such systems 
can be modelled using a variety 
of analytical approximation schemes, including post Newtonian 
expansions, black hole perturbation
theory and the effective one body approach, and this is complemented  
by numerical relativity near the late time coalescence 
where perturbative methods  break down
\cite{Focusissue,Choptuik2015}.

\newpage

\vspace{3pt}

\noindent
{\bf{Problem CMP2:}} 
Determine the predictions of emitted waveforms  for binary black hole 
systems for optimal
detection and parameter extraction.

\vspace{3pt}

Non-vacuum compact binary systems involving at least one neutron
star also produce copious amounts of gravitational waves and are likely 
to lead to intense neutrino and
electromagnetic emission that could also be detected.
However, the simulation of binaries with neutron stars is complicated by 
the need to include nongravitational
physics, and hence analytical techniques are less effective \cite{Choptuik2015}.

Finally, some numerical results have lead to the formulation of new problems in mathematical physics,
some of which have been discussed earlier.
In particular, critical phenomena in gravitational collapse within GR was discovered
numerically \cite{Choptuik2015}.
Families of solutions to the coupled
Einstein-matter equations, labeled by a continuous
parameter $p$, were studied. 
The prescribed initial data depends on  $p$, which 
controls the strength of the (initially imploding) matter in the ensuing gravitational interaction. For a small $p$, gravity is weak during the evolution
and the spacetime remains regular everywhere (for example, in the case of massless radiation,
the radiation will disperse to infinity). For a large $p$, gravity becomes sufficiently strong that
some of the matter is trapped within a singular black hole.
For some critical value of $p$, there is a ``critical'' (self similar) solution corresponding to the threshold of black hole
formation. Evidence to date
suggests that virtually any collapse model that admits black hole formation
will contain such critical behaviour.
Understanding these critical solutions and the ensuing critical behaviour is now an interesting problem within
mathematical GR (especially in the case in which there are no symmetries).

\vspace{3pt}

\noindent
{\bf{Problem CMP3:}}  Understand critical phenomena in gravitational collapse in GR.

\vspace{3pt}

{\em{In summary}}, and in the spirit of {\bf{AL42}} \cite{AL42}, the final 
42 open problems in mathematical physics are: 
H6, S3, S8, S15 and M2, BS1 and BS3, Y117 and Y118,  BS2, BS8 and BS14, 
RB20, RB21, RB32 and RB43, and
the problems P1 -- 17, PF1 -- 3, MPF1 -- 3, and
CMP1 -- 3.

All of these problems are explicitly stated in the text or in the Appendix.
In addition, there are many other open problems referred to in this paper.

%%%%%%%%%%%%%%%%%%%%%%%
\newpage
\section{Appendix: Lists}

\subsection{Hilbert's problems}

The remaining problems are \cite{Hilbert}:

\begin{itemize}

\item H1 The continuum hypothesis (that is, there is no set whose cardinality is strictly between that of the
integers and that of the real numbers). 

\item H2 Prove that the axioms of arithmetic are consistent.  

\item H3 Given any two polyhedra of equal volume, is it always possible to cut the first into finitely many
polyhedral pieces that can be reassembled to yield the second? 

\item H4 Construct all metrics where lines are geodesics.  

\item H5 Are continuous groups automatically differential groups?

\item H7 Is $a^b$ transcendental, for algebraic a $\neq 0,1$ and irrational algebraic b ? 

\item H9 Find the most general law of the reciprocity theorem in any algebraic number field.

\item H10 Find an algorithm to determine whether a given polynomial Diophantine equation with integer coefficients
has an integer solution.  

\item H11 Solving quadratic forms with algebraic numerical coefficients.  

\item H13 Solving 7th degree equations using algebraic functions of two parameters. 

\item H14 Is the ring of invariants of an algebraic group acting on a polynomial ring always finitely generated?

\item H15 Rigorous foundation of Schubert's enumerative calculus.  

\item H17 Express a nonnegative rational function as quotient of sums of squares.  

\item H18 (a) Is there a polyhedron that admits only an anisohedral tiling in three dimensions?  (b) What is the
 densest sphere packing?  

\item H19 Are the solutions of regular problems in the calculus of variations always necessarily analytic?

\item H20 Do all variational problems with certain boundary conditions have solutions? 

\item H21 Proof of the existence of linear differential equations having a prescribed monodromic group.

\item H22 Uniformization of analytic relations by means of automorphic functions. 

\item H23 Further development of the calculus of variations.

\end{itemize}

\subsection{Smale's  problems}

Smale's solved (partially or fully) problems are \cite{Smale}:

\begin{itemize}

\item S6 Finiteness of the number of relative equilibria in celestial mechanics.

\item S7 Distribution of points on the 2-sphere.

\item S11 Is one-dimensional complex-variable dynamics generally hyperbolic? 

\item S12 Centralizers of diffeomorphisms.

\item S14 Lorenz attractor.

\item S17 Solving polynomial equations in polynomial time in the average case.

\end{itemize}

\newpage
\subsection{AL42 problems}

The problems are \cite{AL42}:

\begin{itemize}

\item AL1 Why does conventional physics predict a cosmological constant that is vastly too
large?

\item AL2 What is the dark energy?

\item AL3 How can Einstein gravity be reconciled with quantum mechanics?

\item AL4 What is the origin of the entropy and temperature of black holes?

\item AL5 Is information lost in a black hole?

\item AL6 Did the universe pass through a period of inflation, and if so how and why?

\item AL7 Why does matter still exist?

\item AL8 What is the dark matter?

\item AL9 Why are the particles of ordinary matter copied twice at higher energy?

\item AL10 What is the origin of particle masses, and what kind of masses do neutrinos have?

\item AL11 Does supersymmetry exist, and why are the energies of observed particles so small
compared to the most fundamental (Planck) energy scale?

\item AL12 What is the fundamental grand unified theory of forces, and why?

\item AL13 Are Einstein's GR and standard field theory always valid?

\item AL14 Is our universe stable?

\item AL15 Are quarks always confined inside the particles that they compose?

\item AL16 What are the complete phase diagrams for systems with nontrivial forces, such as
the strong nuclear force?

\item AL17 What new particles remain to be discovered?

\item AL18 What new astrophysical objects are awaiting discovery?

\item AL19 What new forms of superconductivity and superfluidity remain to be discovered?

\item AL20 What further properties remain to be discovered in highly correlated electronic
materials?

\item AL21 What new topological phases remain to be discovered?

\item AL22 What other new phases and forms of matter remain to be discovered?

\item AL23 What is the future of quantum computing, quantum information, and other
applications of entanglement?

\item AL24 What is the future of quantum optics and photonics?

\item AL25 Are there higher dimensions, and if there is an internal space, what is its geometry?

\item AL26 Is there a multiverse?

\item AL27 Are there exotic features in the geometry of spacetime, perhaps including those
which could permit time travel?

\item AL28 How did the universe originate, and what is its fate?

\item AL29 What is the origin of spacetime, why is spacetime four-dimensional, and why is time
different from space?

\item AL30 What explains relativity and Einstein gravity?

\item AL31 Why do all forces have the form of gauge theories?

\item AL32 Why is Nature described by quantum fields?

\item AL33 Is physics mathematically consistent?

\item AL34 What is the connection between the formalism of physics and the reality of human
experience?

\item AL35 What are the ultimate limits to theoretical, computational, experimental, and
observational techniques?

\item AL36 What are the ultimate limits of chemistry, applied physics, and technology?

\item AL37 What is life?

\item AL38 How did life on Earth begin and how did complex life originate?

\item AL39 How abundant is life in the universe, and what is the destiny of life?

\item AL40 How does life solve problems of seemingly impossible complexity?

\item AL41 Can we understand and cure the diseases that afflict life?

\item AL42 What is consciousness?

\end{itemize}

\newpage

\subsection{Simon's problems}

The remaining BS problems are \cite{simon}:

\begin{itemize}

\item BS4 Transport theory:
At some level, the fundamental difficulty of transport theory is that it is a
steady state rather than equilibrium problem, so that the powerful formalism of
equilibrium statistical mechanics is unavailable. A: Fourier's heat law. B: Kubo formula.

\item BS5 Heisenberg models:
Lattice models of statistical mechanics (especially the Ising model)
have been fruitful testing grounds
for ideas in the theory of phase transitions. Four particular questions were postulated, including the  proof
of the Griffiths,
Kelly and Sherman  inequality for classical Heisenberg models.

\item BS6 Existence of ferromagnetism.

\item BS7 Existence of continuum phase transitions.

\item BS9A/B Asymptotic completeness for short range N-body quantum systems
and for Coulomb potentials.

\item BS10 Quantum potential theory:
Basic to atomic and molecular physics is the binding energy of a quantum
mechanical system of electrons interacting with one or more nuclei. Five particular questions were posed.

\item BS11 Existence of crystals: Most materials occur in a crystalline
state at low temperatures. 
Prove the existence of crystals for ensembles of quantum mechanical
atoms (even at zero temperatures) for infinite nuclear masses
with an integer nuclear charge.

\item BS12 Five questions on random and almost periodic potentials.
 
\item BS13 Critical exponent for self-avoiding walks.

\end{itemize}

The two problems in Yau  \cite{Yau1982} referred to earlier are:

\begin{itemize}

\item
Y117 Prove that any YM field on $S^4$ is either self-dual
or antiself-dual.

\item
Y118 Prove that the moduli space of the self-dual fields on $S^4$ 
with a fixed Pontryagin number is connected.

\end{itemize}

\newpage

\subsection{Penrose's problems}

The fourteen unsolved problems in classical
GR by
Penrose (p631 in \cite{Yau1982}) are:

\begin{itemize}

\item 
RP1\&2 Find a suitable quasi-local definition of energy-momentum
in general relativity. And the
more ambitious: Find a suitable quasi-local definition of angular momentum
in general relativity.

\item RP3\&4 Find an ``asymptotically simple''  
(essentially a spacetime in
which every light ray escapes, both in past and future directions, to an
asymptotically flat region) Ricci-flat spacetime which
is not flat, or at least prove that such spacetimes exist.
And the related problem:
Are there restrictions on k for non-stationary
``k-asymptotically simple'' spacetimes, with non-zero mass, which are
vacuum near null infinity.

\item RP5 Find conditions on the Ricci tensor (e.g., satisfying the null convergence 
condition
and the related physically reasonable weak
energy condition)
which ensure that the generators of past and future null infinity (i.e., the 
null geodesic curves lying on
these curves constituting a fibration of null infinity) are infinitely 
long.

\item RP6--8 Assuming appropriate energy conditions hold, show that if a ``cut'' C (a general cross section) of future (or past) null infinity can be spanned by a
spacelike hypersurface,
then the so-called Bondi-Sachs mass defined at C is non-negative.
Does the Bondi-Sachs mass defined on cuts of future null infinity 
have a
well-defined limit as the cuts recede into the past along this limit
agreeing with the mass defined at spacelike infinity?
Show that if the spacetime is assumed not
to be flat everywhere in the region of an appropriate spacelike 
hypersurface, then the
Bondi-Sachs energy-momentum, and also the energy-momentum defined at
spacelike infinity, are future-timelike.

\item RP9\&10 In an asymptotically simple spacetime which is vacuum
near null infinity and for which outgoing radiation is present and falls 
off suitably
near  spacelike and future-timelike infinities, is it necessarily the case 
that spacelike and future-timelike infinities
are nontrivially
related? This leads to: Find a good definition of angular momentum for 
asymptotically
simple spacetimes.

\item RP11 If there is no incoming nor outgoing radiation
and the spacetime manifold is vacuum near future infinity (and, in some 
suitable sense,
near spacelike  infinities) is the manifold necessarily stationary near 
null infinity.

\item RP12 Is cosmic censorship valid in classical
general relativity?

\item RP13 Let S be a spacelike hypersurface  which is compact
with boundary, the boundary consisting of a cut C of future null infinity  together
with a trapped surface (the horizon of the black hole). Then, assuming that the dominant energy condition is satisfied, show that 
there is a lower bound on the ADM mass \cite{ADM} in terms of
the area of S.

\item RP14 Show that there is no vacuum equilibrium configuration
involving more than one black hole.

\end{itemize}

\subsection{Bartnik's problems}

The problems are  \cite{Bartnik}:

\begin{itemize}

\item RB1 Given asymptotically flat initial data and a (future)
trapped surface, prove the existence of (smooth, spherical) apparent horizons.

\item RB2  Is there an analogy between the behaviour of minimal surfaces and the behaviour of
apparent horizons.

\item RB3  Prove that there is an asymptotically flat vacuum initial data set, diffeomorphic to $R^3$, which
contains an apparent horizon.

\item RB4 Determine whether an  asymptotically flat metric on $R^3$ with zero scalar curvature can admit a
minimal 2-sphere (this is the restriction of RB3 to 
time-symmetric initial data).

\item RB5 Give an  explicit example of an apparent horizon that does not persist under the Einstein evolution.

\item RB6 Find conditions on  conformally flat, asymptotically flat metrics with non-negative scalar
curvature which ensure that the manifold has no horizon.

\item RB7  Prove the Penrose inequality.

\item RB8 Determine conditions on the initial data  for a compact
manifold with non-constant mean curvature which ensure the Einstein equation is
solvable.

\item RB9 Describe suitable asymptotic conditions which enable the conformal method
to be used to construct initial data sets on an asymptotically hyperbolic manifold.

\item RB10 Characterise those hyperboloidal initial data arising from a spacetime with
a smooth conformal null infinity.

\item RB11 Classify the various kinds of smoothness properties which hyperboloidal
initial data may have at conformal infinity.

\item RB12 Can the  space of globally hyperbolic, vacuum Einstein metrics on 
$R \times M^3$ have more than one connected component?

\item RB13 What conditions on the stress-energy tensor are needed to show that a
static, asymptotically flat metric is necessarily spherically symmetric?

\item RB14 Find a  local invariant characterisation of the Kerr solution amongst
stationary metrics.

\item RB15 If two disjoint Cauchy surfaces in an asymptotically flat (vacuum)
spacetime are isometric, show that the spacetime is stationary.

\item RB16 Characterize the possible types of singularities which may
occur for solutions to the static and stationary vacuum Einstein equations.

\item RB17 Several questions on the spherically
symmetric Einstein-Yang-Mills equations.

\item RB18 Prove an ``approximate solution'' result for the (vacuum) Einstein
equations in some suitable norm that would 
provide a good way to evaluate approximate/asymptotic and
numerical solutions.

\item RB19 Show that a solution of the linearised (about
Minkowski space) Einstein equations  is close to a (non-flat) exact solution.

\item RB20  Determine the range of validity of the post-Newtonian and post-
Minkowskian asymptotic expansions.

\item RB21 Prove rigorously the existence of a limit in which solutions of the Einstein equations
reduce to Newtonian spacetimes.

\item RB22  Prove the quadrupole radiation formula.

\item RB23 Show that test particles  move on spacetime
geodesics.

\item RB24 In what sense does a Regge spacetime (i.e.,  a piecewise linear 
manifold with piecewise linear  metric \cite{Regge}), and  generally spacetimes constructed by
numerical relativity,  approximate a smooth vacuum
spacetime? 

\item RB25 A problem encountered in numerical relativity is that of ensuring that the
constraint equations are preserved by the evolution \cite{Stewart}.

\item RB26  Prove a uniqueness
theorem for maximal surfaces, assuming only the dominant energy condition.

\item RB27 Rigorously  demonstrate the existence of a constant mean curvature
hypersurface asymptotic to a given cut of future null infinity in an asymptotically flat spacetime.

\item RB28 Show that there is a maximal Cauchy hypersurface of an asymptotically flat spacetime having a Cauchy surface without horizons.

\item RB29 Show that a maximal surface in a ``boost-type domain''  is necessarily
asymptotically flat and must coincide with the maximal slices.

\item RB30 Is there a timelike geodesically complete inextendible Lorentz manifold
satisfying an energy condition and having a partial Cauchy surface which
contains a trapped surface? 

\item RB31 Show that a weak Cauchy surface in a globally hyperbolic spacetime satisfying suitable energy
conditions cannot contain an inextendible null geodesic. 

\item RB32  
Prove that a ``cosmological spacetime'' satisfying the timelike
convergence condition is either timelike geodesically incomplete or it
splits as $R \times M^3$ isometrically (and is thus static).

\item RB33  Prove a singularity theorem assuming the dominant energy condition
rather than the timelike convergence condition.

\item RB34 
Determine the
weakest condition on the  smoothness of the metric in the initial value problem
for maximising  geodesics to have a unique solution.

\item RB35 Prove a long~time existence theorem for the vacuum asymptotically flat Einstein equations
in the maximal slicing gauge \cite{Moncrief81E}.

\item RB36 Determine conditions on asymptotically flat initial data which ensure that the null infinity
of the resulting solution of the initial value problem is sufficiently regular that the Penrose
extended manifold exists \cite{BR6}.

\item RB37 Show the existence of (and construct an exact solution to) the Einstein
vacuum equations with positive mass which has complete smooth  null infinity and  regular timelike infinity.

\item RB38 Find the weakest possible regularity conditions for a metric 
to satisfy the (distributional) vacuum
Einstein equations \cite{HughesKato}.

\item RB39
What are the regularity conditions
for the vacuum Einstein initial value problem  for geodesics which guarantee the existence of a solution, but not
uniqueness.

\item RB40 Demonstrate the  long-time existence of so--called crushing singularities in the constant mean curvature slicing gauge for the
cosmological  vacuum spacetime.

\item RB41 Prove that every globally hyperbolic, maximally extended spacetime
solution of the Einstein or Einstein-Maxwell equations on $R \times S^3$ contains a
maximal hypersurface (and thus also both a big bang and big crunch).
Prove this result for special cases such as the spatially homogeneous Bianchi type IX solutions.

\item RB42 Show that, in an appropriate sense, the set of spacetime metrics which are
smoothly (distributionally?) extendible across compact Cauchy horizons are of ``measure
zero'' in the set of all spacetimes.

\item RB43 Find an exact solution of the Einstein equations which represents two
orbiting bodies. Is the 2-body system unstable in Einstein gravity?

\item RB44 Prove the Belinskii, Khalatnikov and Lifshitz  conjecture.

\item RB45 Show that the only solution of the vacuum (or Einstein-Maxwell) Robinson-Trautman equations
on $S^2 \times R$ with positive mass is the Schwarzschild metric.

\item RB46 Show that a perturbation of the
Schwarzschild (and Kerr) solution decays exponentially (so that the solutions are thus stable).

\item RB47 Show that a cosmological spacetime with constant mean curvature initial data having positive
Ricci 3-curvature has an evolution which preserves Ricci positivity.

\item RB48 Find a sensible notion of quasi-local mass that can be used in non-trivial
black hole theorems.

\item RB49 Show that the set  of asymptotically flat 3-manifolds which satisfy the conditions of the
positive mass theorem has some weak compactness property (and what regularity might be
expected in the limit manifold?)

\item RB50  Prove the static metric
conjecture \cite{BR6}.

\item RB51 Construct a general proof of the positive mass theorem that does not require the existence of a foliation with special properties.

\item RB52  Show that the Bartnik quasi-local mass \cite{BR6} is strictly positive for non-flat
data and that the Penrose quasi-local mass is non-negative for reasonable
data.

\item RB53 Explain the relation between the various definitions of quasi-local mass.

\end{itemize}

\newpage

\subsection{Lists of lists}

\begin{itemize}

\item List of acronyms used in this paper.

\item Top 10 movies.

\item Top 10 songs.

\item Further list of lists$^1$.

\end{itemize}

$^1$ A list of the top 10 books, perhaps subdivided into popular books, popular science books and technical science books in mathematical physics, has proven more problematic to formulate.

\subsubsection{Acronyms}

Let us present a  list of commonly used acronyms:
\begin{itemize}

\item AdS:
Anti-de Sitter spacetime. 

\item BKL: Belinskii, Khalatnikov and Lifshitz. 

\item CFT: quantum field theory.

\item CMB: cosmic microwave background.

\item EYM: Einstein-Yang-Mills. 

\item FLRW:
Friedmann--Lemaitre--Robertson--Walker.  

\item GR: general relativity.

\item LQC: loop quantum cosmology. 

\item LQG: loop quantum gravity. 

\item QFT: quantum field theory.

\item QG: quantum gravity.

\item YM: Yang-Mills. 

\item 4D: four dimensions.

\end{itemize}

The abbreviations for the various lists are: 
H (Hilbert), S (Smale), M (Millennium), 
AL (Allen and Lidstrom), BS (Simon),
Y (Yau), RB (Bartnik), RP  (Penrose), Ph (theoretical physics),
PF (personal favorites),  P (open mathematical physics problems in the contemporary 
fields of GR, the quantum realm and in cosmology), MPF (mathematical physics problems in the future),
and CMP (computationally motivated problems). All of the problems referred to are
explicitly stated in the text or in the Appendix.

%%%Other useful information: SH

\newpage

\subsubsection{Top 10 movies}

\begin{itemize}

\item 
Ramanujan (2014).
Director: Gnana Rajasekaran.

\item Einstein and Eddington (2008). 
Director: Philip Martin.

\item Infinity (1996).
Director: Matthew Broderick.

\item  The Theory of Everything (2014). 
Director: James Marsh

\item A Beautiful Mind (2001).
Director: Ron Howard.

\item Proof (2005).
Director: John Madden.

\item Good Will Hunting (1997). 
Director: Gus Van Sant.

\item
Pi (1998).
Director: Darren Aronofsky.

\item The Imitation Game (2014).
Director: Morten Tyldum.

\item
Fermat's Room
(2007).  Directors: Luis Piedrahita and Rodrigo Sopena.

\end{itemize}

\subsubsection{Top 10 songs}

\begin{itemize}

\item Eric Idle (Monty Python) -- Galaxy Song.

\item The Bare Naked Ladies -- The Big Bang.

\item
Kate Bush -- Pi.

\item
Jack Black -- Math Song.

\item
Jarvis Cocker -- Quantum Theory.

\item
Nick Cave and The Bad Seeds -- Higgs Boson Blues.

\item
They Might Be Giants -- Why Does the Sun Shine? 

\item 
MC Hawking -- $E=mc^2$.

\item
One Direction -- Maths Song.

\item
Bjork -- Mutual Core.

\item
Louis Armstrong -- What a Wonderful World $^2$.

\end{itemize}

All of these songs are available on You Tube.

$^2$
No list of top songs should ever exclude this song.

\newpage

\section*{Acknowledgements}  
I would like to thank Lars Andersson, Robert van den Hoogen and Claes Uggla
for a detailed reading of an earlier version of the manuscript, and  Tim Clifton, Luis Lehner and Frans Pretorius for helpful comments. Financial support was provided by  NSERC of Canada.

\end{document}